\documentclass[11pt]{article}
\usepackage{amsmath,amsfonts,amscd,bezier}
\usepackage{latexsym,amssymb}
\usepackage[latin1]{inputenc}
\usepackage[english]{babel}

\newtheorem{theorem}{Theorem}[section]
\newtheorem{lemma}[theorem]{Lemma}
\newtheorem{proposition}[theorem]{Proposition}
\newtheorem{corollary}[theorem]{Corollary}
\newtheorem{definition}[theorem]{Definition}
\newtheorem{example}[theorem]{Example}
\newtheorem{remark}[theorem]{Remark}

\newcommand{\Dem}{\noindent{\it Proof.}\ \ }
\newcommand{\cqd}{{\hfill $\rule{2mm}{2mm}$}\vspace{.5cm}}

\def\Re{\mathbb{R}}

\begin{document}

\title{\bf Normal forms of $\omega$-Hamiltonian vector fields with symmetries}

\author{Patr\'{\i}cia~H.~Baptistelli, Maria Elenice R. Hernandes\\
{\small Departamento de Matem\'atica, Universidade de Maring\'a}\\
{\small Av. Colombo, 5790, 87020-900, Maring\'a - PR, Brazil \footnote{Email address: phbaptistelli@uem.br (corresponding author), merhernandes@uem.br}}
\and Eralcilene Moreira Terezio\\
{\small Instituto Latino-Americano de Ci\^encias da Vida e da Natureza}\\ 
{\small Universidade Federal da Integra\c{c}\~ao Latino-Americana} \\ {\small Av. Tancredo Neves, 6731,  85867-900, Foz do Iguaçu - PR, Brazil \footnote{Email address: eralcilene.terezio@unila.edu.br}}}

\date{}
\maketitle

\begin{abstract}
In this paper, we present algebraic tools to obtain normal forms of $\omega$-Hamiltonian vector fields under a semisymplectic action of a Lie group, by taking into account the symmetries and reversing symmetries of the vector field. The normal forms resulting from the process preserve the Hamiltonian condition and the types of symmetries of the original vector field. Our techniques combine the classical method of normal forms of Hamiltonian
vector fields with the invariant theory of groups.
\end{abstract}


\noindent {2020 AMS Classification: 37C79, 37J40, 58D19} \\

\section{Introduction}
\label{intro}
There are several methods to understand the local qualitative behavior of a vector field and one of them is the normal form theory. Some authors, as in \cite{buzzi2,mauricio,Ricardo2013,Ri,Me}, use this theory to study limit cycles, families of periodic orbits, relative equilibria, relative periodic solutions, etc.

The classic way to obtain a normal form of a vector field consists of performing changes of coordinates around a singular point (namely $x=0$) on a formal power series of the vector field. More precisely, we consider perturbations of the identity of the form $x = y+\varphi^k(y)$, with $k\geq2$ and $\varphi^k$ a homogeneous polynomial of degree $k$, which preserve the linear part and allow to eliminate terms of degree $k\geq 2$ of the  original vector field, by resulting in a formally conjugate vector field written in a more convenient way. For each $k\geq2$, the truncated equation of degree $k$ of the associated power series is called a normal form of order $k$ of the vector field.

Over the years, some methods have been developed to provide normal forms as described above. Belitskii \cite{Be,Be2} reduces this problem to the computation of solutions of a partial differential equation. Elphick {et al.} \cite{elphic} provide an algebraic method, in which it is possible to obtain a normal form by choosing non-linear terms that commute with the action of the one-parameter group\begin{equation}\label{grupoS}
	\mathbf{S}=\overline{\left\{e^{sL^T}\hspace{-0.23cm}: \ s\in{\mathbb{R}}\right\}},
\end{equation}where the bar denotes the closure of the set and $L$ is the linearization of the vector field at the origin. In this case, the normal form process introduces additional symmetries into the problem since every vector field is formally conjugated to a vector field with symmetry group $\mathbf{S}$. Several works deal with vector fields  that present symmetric geometric configurations, as we can see in \cite{Bap,Bap2,broer,Golu2,Montaldi}. In particular, in the same way as Elphick et al. \cite{elphic} and based on algebraic invariant theory, Baptistelli, Manoel and Zeli  present a method (\cite[Theorem 4.7]{Bap}) to obtain normal forms of a reversible equivariant vector field that preserve the symmetries of the original vector field. 

The class of Hamiltonian vector fields, and consequently Hamiltonian systems, deserves special attention due to its historical contribution in Mathematics and its applications in Physics, among others. Given $(V,\omega)$ a symplectic real vector space of dimension $2n$ and $H:V\to\mathbb{R}$ a smooth function, the $\omega$-Hamiltonian vector field associated with $H$ can be written as $X(x) = ([\omega]^{-1})^{T}\nabla H(x)$ for all $x\in V$, where $[\omega]$ denotes the matrix of $\omega$ relative to a basis of $V$ (see Definition \ref{deficampoomegahamil}). In symmetric Hamiltonian context, several works have contributed in different aspects, such as \cite{Alomair, buono, buzzi1, Montaldi}. In \cite{Montaldi}, Montaldi, Roberts and Stewart deal with Hamiltonian systems on $(\mathbb{R}^{2n},\omega_0)$, where $\omega_0$ is the canonical bilinear form, under an orthogonal action of a Lie group that admits an involutory reversible symmetry. In \cite{buzzi1}, Buzzi and Teixeira also analyze the dynamics of reversible Hamiltonian vector fields on $(\mathbb{R}^{2n},\omega_0)$  with respect to a symplectic involution. More recently, Alomair and Montaldi \cite{Alomair} study the existence of families of periodic solutions in two classes of Hamiltonian systems with involutory symmetry (of types SR and AE according to Table \ref{tabelatipos}). In that paper, the authors consider two group homomorphisms and classify the symmetries into up four types. 

The aim of this work is presenting an algebraic method to obtain normal forms of $\omega$-Hamiltonian vector fields on $(V,\omega)$ under a semisymplectic action of a Lie group $\Gamma$ (see Definition \ref{definiacaosemi}), by taking into account the symmetries and reversing symmetries of $X$. For this purpose, we assume the existence of two group homomorphisms $\sigma_1, \sigma_2: \Gamma \to \mathbb{Z}_2$, where $\sigma_1$ is related to semisymplectic action of $\Gamma$ and $\sigma_2$ establishes the symmetries of the vector field. In this context, we consider that the matrix of the symplectic form $\omega$ may possibly be different from the matrix of the canonical symplectic form given by \begin{equation}\label{Jnporn}
	J = \left[\begin{array}{cc}
		0&I_n\\-I_n&0
	\end{array}\right],
\end{equation} where $I_n$ is the identity matrix of order $n$. One of the reasons for this is that the symmetries of vector fields may not be preserved by symplectomorphisms. In addition, some works (as \cite{Robinson,Sev}) also deal with Hamiltonian vector fields for matrices of the symplectic form different from $J$.

The method we present gives us a normal form that preserves the types of symmetries and Hamiltonian condition of the original vector field. In this case, there is a computational advantage compared to the method presented in \cite{Bap}, since a normal form of a Hamiltonian vector field can be obtained by means of a normal form of the associated Hamiltonian function (see, for instance, \cite{Chow}).

This paper is organized as follow: in Section \ref{sec:2preliminar} we introduce some notations and basic notions of Hamiltonian and equivariant contexts. In Section \ref{sec:semisymplectic} we use the semisymplectic group to recognize a semisymplectic action on $(V,\omega)$ (Theorem \ref{acao subgrupo de Omega}). In Section \ref{sec:3simetrias}, we describe the types of symmetry in Hamiltonian context with respect to a (Lie) group of symmetries $\Gamma$ and its action on $(V,\omega)$.  In Section \ref{secaoformasnormaisNOVA} we present the classic normal form theory of $\omega$-Hamiltonian vector fields and give a method to compute these normal forms (Subsection \ref{secaoDK}). The Section \ref{secaoformasnormaissimetrias} deals with the main issue of the paper, by presenting an algebraic method to obtain normal forms of $\omega$-Hamiltonian vector fields with symmetries (Theorem \ref{fnormalnossa} and Corollary \ref{fnormalnossa2}). We illustrate such a method with two examples (Subsections \ref{sec:example} and \ref{sec:example2}) and we finish the paper with a section of some technical proofs (Section \ref{adsubsecaopropriedades}).

\section{Preliminary concepts}
\label{sec:2preliminar}

\subsection{Preliminaries on symplectic geometry}

A symplectic vector space is a pair $(V,\omega)$, where $V$ is a finite-dimensional real vector space and $\omega:V\times V\to\mathbb{R}$ is an alternating and non-degenerate bilinear form. It is possible to prove that the dimension of $V$ is even and the matrix $[\omega]$ of $\omega$ relative to a basis of $V$ is a skew-symmetric and invertible matrix. The next definition presents our main object of study and it can naturally be extended to a symplectic manifold.

\begin{definition}\label{deficampoomegahamil}
	Let $(V,\omega)$ be a symplectic vector space and $H:V\to\mathbb{R}$ a smooth function. The $\omega$-Hamiltonian vector field associated with $H$ is the unique vector field $X_H:V\to V$ such that $\omega(X_H(x),\cdot) = dH_x$, for all $x\in V$. In other words, $X_H$ is written as
	\begin{equation}\nonumber
		X_H(x) = ([\omega]^{-1})^{T}\nabla H(x)
	\end{equation}for all $x\in V$, where $T$ denotes the transpose of the matrix and $\nabla H$ is the gradient vector of $H$. The function $H$ is called a Hamiltonian function.
\end{definition}

We intend performing changes of coordinates in an $\omega$-Hamiltonian vector field by preserving its Hamiltonian structure. Changes of coordinates called symplectic or antisymplectic satisfy this property (see \cite[Theorem 4.2]{BRHT2021}).

\begin{definition}\label{definicaolambdasimpletica} Let $(V,\omega)$ be a symplectic vector space. A differentiable map $\xi:V\to V$ is called symplectic if $\xi^*\omega = \omega$ and it is called antisymplectic if 
	$\xi^*\omega = -\omega$, where $\xi^*\omega$ is the pullback of $\omega$ by $\xi$.
\end{definition}

If $(V,\omega)$ has dimension $2n$, then the matrices $B\in\mathbb{M}_{2n}(\mathbb{R})$ associated with symplectic linear maps satisfy the equality $B^T[\omega] B = [\omega]$ and are called $\omega$-symplectic matrices (see \cite[Proposition 2.3]{BRHT2021}). The set of all $\omega$-symplectic matrices is a Lie group denoted by $Sp_{\omega}(n;\mathbb{R})$ and called $\omega$-symplectic group. In addition, matrices $B$ which are associated with antisymplectic linear maps are called $\omega$-antisymplectic matrices and satisfy $B^T[\omega] B = -[\omega]$. The set of all $\omega$-antisymplectic matrices is  denoted by $Sp_{\omega}^{-1}(n;\mathbb{R})$ and it does not have a group structure.

\begin{remark}\label{deltaSpnA}
	We highlight the following properties of the set $Sp_{\omega}^{-1}(n;\mathbb{R})$: if $B\in Sp_{\omega}^{-1}(n;\mathbb{R})$ then $B^{-1}\in Sp_{\omega}^{-1}(n;\mathbb{R})$; the product of an even number of elements in $Sp_{\omega}^{-1}(n;\mathbb{R})$ belongs to $Sp_{\omega}(n;\mathbb{R})$, while the product of an odd number of elements in $Sp_{\omega}^{-1}(n;\mathbb{R})$ still belongs to $Sp_{\omega}^{-1}(n;\mathbb{R})$.
\end{remark}

In \cite{BRHT2021} we present a study of the symplectic algebra for an arbitrary symplectic vector space $(V,\omega)$, by introducing the concept of $\omega$-symplectic group and characterizing its Lie algebra as a useful tool in the recognition of $\omega$-Hamiltonian vector fields. More precisely, we define the $\omega$-semisymplectic group by the disjoint union
\begin{equation}
	\label{eq:omegan}
	\Omega_n = Sp_{\omega}(n;\mathbb{R}) \ \dot{\cup} \  Sp_{\omega}^{-1}(n;\mathbb{R}).\end{equation} The Lie algebra of $\Omega_n$, which coincides with the Lie algebra of $Sp_{\omega}(n;\mathbb{R})$, is given by the set

\begin{equation}
\label{Liealgebra}
{\mathfrak{sp}}_{\omega}(n;{\mathbb{R}}) = \left\{L\in {\mathbb{M}}_{2n}({\mathbb{R}}): L^T[\omega]+[\omega]L = 0 \right\}.
\end{equation}

\noindent The elements of this algebra are called $\omega$-Hamiltonian matrices.

%
%
%

\subsection{Preliminaries on invariant theory}


Consider  $V$ a finite-dimensional real vector space and $\Gamma$ a Lie group acting on $V$ by means of a differentiable action $(\gamma, x) \mapsto \gamma x$, for all $\gamma \in \Gamma$ and $x \in V$. Consider ${\rho}_{\gamma}:V\to V$ the diffeomorphism induced by this action, given by ${\rho}_{\gamma}(x)=\gamma x$. If $X:V\to V$ is a vector field, for each $\gamma\in\Gamma$ we denote by $(\rho_{\gamma})_{*}X: V \to V$ the pushforward of $X$ by ${\rho}_{\gamma}$, namely
$(\rho_{\gamma})_{\ast}X(y) = d(\rho_{\gamma})_x(X(x))$ for $x = \rho_{\gamma}^{-1}(y)$. We say that $X$ is $\Gamma_{\sigma}$-equivariant if there exists a group homomorphism $\sigma:\Gamma\to\mathbb{Z}_2 = \{\pm1\}$ satisfying $(\rho_{\gamma})_{*}X = \sigma(\gamma)X,$ that is,
$$X(\rho_{\gamma}(x)) = \sigma(\gamma)d(\rho_{\gamma})_x(X(x))$$
for all $\gamma\in\Gamma$ and $x \in V$. When $\sigma$ is trivial, $X$ is called $\Gamma$-equivariant. An element $\gamma\in\Gamma$ is a symmetry of $X$ if $\gamma\in\ker\sigma$ and it is a reversible symmetry of $X$ if $\gamma\in\Gamma\backslash\ker\sigma$. 

A function $f:V\to\mathbb{R}$ is called $\Gamma_{\sigma}$-invariant if
	\[f(\rho_{\gamma}(x))=\sigma(\gamma)f(x),\]
for all $\gamma\in\Gamma$ and $x\in
V$. If $\sigma$ is trivial, then $f$ is called $\Gamma$-invariant.


Denote by $\mathcal{P}$ and $\overrightarrow{\mathcal{P}}$ the vector spaces of all polynomial functions $V\to\mathbb{R}$ and polynomial maps $V\to V$, respectively. We introduce the following notations:\begin{equation}\label{conjSIMETRIAS}
	\begin{array}{ccl}
		{{\mathcal{P}}}(\Gamma)&= &\{f\in\mathcal{P}: f(\rho_{\gamma}(x)) = f(x), \ \forall\gamma\in\Gamma, \ x\in V\};\\
		
		{{\mathcal{P}}_{\sigma}}(\Gamma) &= &\{f\in\mathcal{P}: f(\rho_{\gamma}(x)) =\sigma(\gamma) f(x), \ \forall\gamma\in\Gamma, \ x\in V\};\\
		
		\overrightarrow{\mathcal{P}}(\Gamma) &= &\{F\in\overrightarrow{\mathcal{P}}: F(\rho_{\gamma}(x)) =\rho_{\gamma}(F(x)), \ \forall\gamma\in\Gamma, \ x\in V\};\\
		
		{\overrightarrow{\mathcal{P}}_{\sigma}}(\Gamma)& = &\{F\in\overrightarrow{\mathcal{P}}: F(\rho_{\gamma}(x)) =\sigma(\gamma) \rho_{\gamma}(F(x)), \ \forall\gamma\in\Gamma, \ x\in V\}.
	\end{array}
\end{equation}
The set ${{\mathcal{P}}}(\Gamma)$ has a ring structure and ${{\mathcal{P}}_{\sigma}}(\Gamma)$ is a module over ${{\cal{P}}}(\Gamma)$. If $\Gamma$ acts linearly on $V$, the sets ${\overrightarrow{\mathcal{P}}}(\Gamma)$ and ${\overrightarrow{\mathcal{P}}_{\sigma}}(\Gamma)$ are also modules over ${{\cal{P}}}(\Gamma)$.

For our purposes, given two groups $\Gamma_1$ and $\Gamma_2$, an action of the semidirect product $\Gamma_1\rtimes\Gamma_2$ on $V$ can be given by the mapping $(\Gamma_1\rtimes\Gamma_2)\times V\to V$ defined by $(\gamma_1,\gamma_2)x=\gamma_1(\gamma_2 x)$. For $j=1$,$2$, if $\beta_j:\Gamma_j\to\mathbb{Z}_2$ is a group homomorphism, we define the homomorphism ${\beta}:\Gamma_1\rtimes\Gamma_2\to\mathbb{Z}_2$ by \begin{equation}\label{6barra}{\beta}(\gamma_1,\gamma_2)=\beta_1(\gamma_1)\beta_2(\gamma_2).
\end{equation}

The following result relates the invariant theory for $\Gamma_1\rtimes\Gamma_2$ with the invariant theory for each group, $\Gamma_1$ and $\Gamma_2$. This is a natural extension of \cite[Proposition 3.2]{Bap}, where the authors consider $\beta_1$ in (\ref{6barra}) as the trivial homomorphism.

\begin{proposition}\label{propotrintadois}Let $\Gamma_1$ and $\Gamma_2$ be Lie groups acting linearly on $V$ and ${\beta}$ defined as in (\ref{6barra}). Then
	\begin{enumerate}
		\item ${{\mathcal{P}}}(\Gamma_1\rtimes\Gamma_2) = {{\mathcal{P}}}(\Gamma_1) \cap {{\mathcal{P}}}(\Gamma_2)$;
		\item ${\overrightarrow{\mathcal{P}}}(\Gamma_1\rtimes\Gamma_2) = {\overrightarrow{\mathcal{P}}}(\Gamma_1) \cap {\overrightarrow{\mathcal{P}}}(\Gamma_2)$;
		\item ${{\mathcal{P}}_{{\beta}}}(\Gamma_1\rtimes{\Gamma_2}) = {{\mathcal{P}}_{\beta_1}}({\Gamma}_1) \cap {{\mathcal{P}}_{\beta_2}}({\Gamma}_2)$;
		\item ${\overrightarrow{\mathcal{P}}_{{\beta}}}(\Gamma_1\rtimes\Gamma_2) = {\overrightarrow{\mathcal{P}}_{\beta_1}}({\Gamma}_1) \cap {\overrightarrow{\mathcal{P}}_{\beta_2}}({\Gamma}_2)$.
	\end{enumerate}
\end{proposition}

\section{Semisymplectic actions}\label{sec:semisymplectic}

Semisymplectic actions on a finite-dimensional symplectic vector space are the natural actions in the study of Hamiltonian systems with symmetries and reversing symmetries (see, for instance, \cite{Alomair, Montaldi, Claudia}). For what follows, $\Gamma$ is a Lie group acting on a $2n$-dimensional symplectic vector space $(V,\omega)$. 



\begin{definition}\label{definiacaosemi}
	An action of $\Gamma$ on $(V,\omega)$ is called {$\sigma$-semisymplectic} if there exists a group homomorphism $\sigma:\Gamma\to {\mathbb{Z}}_2$ such that the pullback $(\rho_{\gamma})^*\omega$ satisfies $(\rho_{\gamma})^*\omega = \sigma(\gamma)\omega$, that is, \begin{equation}\nonumber
		\omega(d(\rho_{\gamma})_x\mathbf{u},d(\rho_{\gamma})_x\mathbf{v}) = \sigma(\gamma)\omega(\mathbf{u},\mathbf{v}),
	\end{equation}for all $\gamma\in\Gamma$ and $x,\mathbf{u},\mathbf{v}\in V$. If $\sigma$ is the trivial homomorphism, we say that the action of $\Gamma$ on $V$ is symplectic or $\Gamma$ {acts symplectically} on $V$.
\end{definition}

An element $\gamma\in\ker\sigma$ is called {symplectic} and $\gamma\in\Gamma\setminus\ker\sigma$ is called {antisymplectic}. These nomenclatures are suitable, since $\gamma\in\ker\sigma$ if and only if $\rho_{\gamma}$  is a symplectic map and $\gamma\in\Gamma\setminus\ker\sigma$ if and only if $\rho_{\gamma}$ is an antisymplectic map. 

If $\sigma$ is surjective, then the quotient group $\Gamma/ \ker\sigma$ is isomorphic to $\mathbb{Z}_2$. Thus, given $\delta\in\Gamma\setminus\ker\sigma$, we have $\Gamma\setminus\ker\sigma = \delta\ker\sigma$ as the unique non-trivial left-coset of $\ker\sigma$. Hence\[\Gamma = \ker\sigma \ \dot{\cup} \ \delta\ker\sigma.\] 


If the action of $\Gamma$ on $V$ is linear, the next result characterizes subgroups of the $\omega$-symplectic and $\omega$-semisymplectic groups in terms of properties of this action. For this purpose, we identify the group  $GL(V)$ of invertible linear operators on $V$ with the group $GL(2n)$ of invertible matrices of order $2n$.

\begin{theorem}\label{acao subgrupo de Omega}
	Let $\Gamma$ be a Lie group acting linearly on a $2n$-dimensional symplectic vector space $(V,\omega)$ and $\rho:\Gamma\to GL(2n)$ the representation of $\Gamma$ on $V$. Then $\rho(\Gamma)$ is a subgroup of $\Omega_n$ defined in (\ref{eq:omegan}) if and only if the action of $\Gamma$ on $(V,\omega)$ is $\sigma$-semisymplectic, for some group homomorphism $\sigma:\Gamma\to\mathbb{Z}_2$. In particular, $\rho(\Gamma)$ is a subgroup of $Sp_{\omega}(n;{\mathbb{R}})$ if and only if the action of $\Gamma$ on $(V,\omega)$ is symplectic. 
\end{theorem}
\Dem Suppose $\rho(\Gamma)$ is a subgroup of $\Omega_n$ and identify $\rho(\gamma) :=\rho_{\gamma}$ with its matrix, for each $\gamma\in\Gamma$. Then $\rho_{\gamma}\in Sp_{\omega}(n;{\mathbb{R}})$ or $\rho_{\gamma}\in Sp_{\omega}^{-1}(n;{\mathbb{R}})$. 
	
	By \cite[Proposition 2.3]{BRHT2021}, if $\rho(\Gamma)\subset Sp_{\omega}(n;{\mathbb{R}})$, then $\rho_{\gamma}$ is a symplectic linear operator for all $\gamma\in\Gamma$, that is, $(\rho_{\gamma})^*\omega = \omega$. Thus, the action of $\Gamma$ on $(V,\omega)$ is symplectic. 
	
	If $\rho(\Gamma)\cap Sp_{\omega}^{-1}(n;{\mathbb{R}})\neq\emptyset$, we define the map $\sigma:\Gamma\to\mathbb{Z}_2$ by\[\sigma(\gamma) = \left\{\begin{array}{rl}1, &\text{if} \ \rho_{\gamma}\in Sp_{\omega}(n;{\mathbb{R}})\\
		-1, &\text{if} \ \rho_{\gamma}\in Sp_{\omega}^{-1}(n;{\mathbb{R}})
	\end{array}\right..\]Note that $\sigma$ is a group epimorphism (Remark \ref{deltaSpnA}). Furthermore, $\rho_{\gamma}$ is a symplectic map if $\rho_\gamma\in Sp_{\omega}(n;{\mathbb{R}})$, whence $(\rho_{\gamma})^*\omega = \omega=\sigma(\gamma)\omega$. Similarly, it follows from \cite[Proposition 2.3]{BRHT2021} that  $\rho_{\gamma}$ is an antisymplectic map if $\rho_\gamma\in Sp_{\omega}^{-1}(n;{\mathbb{R}})$, whence $(\rho_{\gamma})^*\omega = -\omega = \sigma(\gamma)\omega$.
	
	Conversely, if $\Gamma$ is a Lie group acting linearly and $\sigma$-semisymplectically on $(V,\omega)$, for some homomorphism $\sigma:\Gamma\to\mathbb{Z}_2$, then $\rho_{\gamma}$ is a symplectic map if $\gamma\in\ker\sigma$ and antisymplectic if $\gamma\in\Gamma\setminus\ker\sigma$. If $\gamma\in\ker\sigma$ then $\rho_{\gamma}\in Sp_{\omega}(n;{\mathbb{R}})$ and if $\gamma\in\Gamma\setminus\ker\sigma$ then $\rho_{\gamma}\in Sp_{\omega}^{-1}(n;{\mathbb{R}})$. Hence, $\rho(\Gamma)$ is a subgroup of $\Omega_n$ such that $\rho(\ker\sigma)\subset Sp_{\omega}(n;{\mathbb{R}})$. In particular, if the action of $\Gamma$ on $(V,\omega)$ is symplectic, we can consider $\sigma$ as the trivial homomorphism and, therefore, $\rho_{\gamma}\in Sp_{\omega}(n;{\mathbb{R}})$ for all $\gamma\in\Gamma$. \cqd 

As we mentioned, an important group in the study of dynamical systems which presents symmetric geometric configurations is the linear Lie group $\mathbf{S}$ given in (\ref{grupoS}). When $[\omega] = J$ as in (\ref{Jnporn}) and $L\in{\mathfrak{sp}}_{\omega}(n;{\mathbb{R}})$, it is known that $\mathbf{S}$ is a subgroup of $Sp_{\omega}(n;\mathbb{R}).$ In particular, the authors in \cite{Montaldi} use this property to obtain some results for Hamiltonian vector fields. In a more general context, we characterize when $\mathbf{S}$ is a subgroup of $Sp_{\omega}(n;\mathbb{R})$ in terms of the matrices in ${\mathfrak{sp}}_{\omega}(n;{\mathbb{R}})$.

\begin{proposition}\label{exacaosimpleticasoS2}{Let $(V,\omega)$ be a symplectic vector space and consider the group $\mathbf{S}$ defined in (\ref{grupoS}). Then $\mathbf{S}$ is a subgroup of $Sp_{\omega}(n;\mathbb{R})$ if and only if $L^T\in{\mathfrak{sp}}_{\omega}(n;{\mathbb{R}})$.
}\end{proposition}
\Dem We can consider the Lie algebra of $Sp_{\omega}(n;\mathbb{R})$ as\begin{equation}\label{eqtrianguloElenice}
		{\mathfrak{sp}}_{\omega}(n;{\mathbb{R}}) = \left\{L\in {\mathbb{M}}_{2n}({\mathbb{R}}): e^{tL}\in Sp_{\omega}(n;\mathbb{R}), \ \forall t\in{\mathbb{R}}\right\},
	\end{equation} since $Sp_{\omega}(n;\mathbb{R})$ is a closed subgroup of $GL(2n)$. If $\mathbf{S} \subset Sp_{\omega}(n;\mathbb{R})$, then $e^{sL^T} \in Sp_{\omega}(n;\mathbb{R})$, for all $s\in\mathbb{R}$. Thus $L^T\in{\mathfrak{sp}}_{\omega}(n;{\mathbb{R}})$ by (\ref{eqtrianguloElenice}). On the other hand, if $L^T\in{\mathfrak{sp}}_{\omega}(n;{\mathbb{R}})$, then $e^{sL^T} \in Sp_{\omega}(n;\mathbb{R})$, for all $s\in\mathbb{R}$, that is, $\{e^{sL^T}: s\in\mathbb{R}\} \subset Sp_{\omega}(n;\mathbb{R})$. We claim that $\mathbf{S}$ is also a subgroup of $Sp_{\omega}(n;\mathbb{R})$. Indeed, given $B\in{\bf S}$, there exists a sequence $(B_m)\subset
	\{ e^{sL^T}: s\in\mathbb{R}\}$ such that $\lim B_m = B$. Then\[B^T[\omega]B = (\lim B_m)^T[\omega](\lim B_m) =  \lim \left(B_m^T[\omega]B_m \right)  = \lim [\omega] = [\omega],\]that is, $B\in Sp_{\omega}(n;\mathbb{R})$. 	\cqd  

\begin{corollary}\label{exacaosimpleticasoS}
	Let $(V,\omega)$ be a $2n$-dimensional symplectic vector space such that $[\omega]^2 = -I_{2n}$. If $L\in {\mathfrak{sp}}_{\omega}(n;{\mathbb{R}})$, then $\mathbf{S}$ is a subgroup of $Sp_{\omega}(n;\mathbb{R})$.
\end{corollary}

\Dem Note that if $[\omega]^2 = -I_{2n}$, then by (\ref{Liealgebra}) we have that ${\mathfrak{sp}}_{\omega}(n;{\mathbb{R}}) = \{L\in{\mathbb{M}}_{2n}({\mathbb{R}}): [\omega] L^T[\omega]  = L\}$. In this case, if $L\in{\mathfrak{sp}}_{\omega}(n;{\mathbb{R}})$, then $L^T\in{\mathfrak{sp}}_{\omega}(n;{\mathbb{R}})$. Thus, the result is a consequence of the previous proposition. \cqd  

\begin{example}\label{Snaosimpleticoex}{\rm Let $(\mathbb{R}^4,\omega)$ be the symplectic vector space and $L$ in the definition of $\mathbf{S}$ such that \begin{equation}\nonumber
			[\omega]=\left[\begin{array}{cccc}
				0&1&0&2\\
				-1&0&-1&0\\
				0&1&0&1\\
				-2&0&-1&0
			\end{array}\right] \quad  \textrm{ and} \quad L = \left[\begin{array}{cccc}
				-1&1&-1&2\\
				3&0&4&1\\
				-1&2&0&2\\
				3&1&1&1
			\end{array}\right].
		\end{equation}Since $L^T\notin{\mathfrak{sp}_{\omega}(2,\mathbb{R})}$, Proposition \ref{exacaosimpleticasoS2} ensures that $\mathbf{S}$ is not a subgroup of $Sp_{\omega}(2;\mathbb{R})$.}\end{example}

\section{Types of symmetry on Hamiltonian context}\label{sec:3simetrias}

Consider $\sigma_1,\sigma_2:\Gamma\to {\mathbb{Z}}_2$ two group homomorphisms. Assume that $\Gamma$ acts $\sigma_1$-semisympletically on $(V,\omega)$. In this section, we relate symmetries of a $\Gamma_{\sigma_2}$-equivariant $\omega$-Hamiltonian vector field with symmetries of its associated Hamiltonian function.

The  motivation for this correspondence are some works in the literature. For instance, in \cite[Proposition 3]{buzzi1} Buzzi and Teixeira consider reversible Hamiltonian vector fields under the action of a linear involution $\delta$ and verify that their associated Hamiltonian functions $H$ are $(\mathbb{Z}_2)_{\sigma}$-invariant, with $\sigma(\delta) = -1$. Alomair and Montaldi \cite{Alomair} extend this result by considering an orthogonal action of a matrix group $\Gamma$ on $({\mathbb{R}}^{2n},\omega_0)$ and show that there are up to four types of symmetry for Hamiltonian vector fields. Theorem \ref{teogenebuzzi1} generalizes these approach, in the sense that we impose neither symplectic coordinates, nor linearity and orthogonality of the action of $\Gamma$ on $(V,\omega)$.

From an algebraic point of view, the four types of symmetry $\gamma\in\Gamma$ as presented in  \cite{Alomair} are: symplectic equivariant (SE), when $\sigma_1(\gamma) = \sigma_2(\gamma) = 1$; symplectic reversible (SR), when $\sigma_1(\gamma) = -\sigma_2(\gamma)=1$; antisymplectic equivariant (AE), when $\sigma_1(\gamma) = -\sigma_2(\gamma)=-1$ and antisymplectic reversible (AR), when $\sigma_1(\gamma) = \sigma_2(\gamma)=-1$. 

In the next result we denote by $[G:H]$ the index of the subgroup $H$ of $G$ in $G$.
\begin{proposition}\label{proptipos}Let $\Gamma$ be a Lie group acting on a symplectic vector space $(V,\omega)$ and $\sigma_1,\sigma_2:\Gamma\to{\mathbb{Z}}_2$ two distinct epimorphisms. If the action of $\Gamma$ on $(V,\omega)$ is $\sigma_1$-semisymplectic and $X$ is an $\Gamma_{\sigma_2}$-equivariant $\omega$-Hamiltonian vector field, then $X$ presents all types of symmetry: SE, SR, AE and AR.
\end{proposition}
\Dem Clearly the identity of $\Gamma$ is SE. Since $\sigma_1$ and $\sigma_2$ are distinct we have $\ker\sigma_1\neq\ker\sigma_2$. Moreover $\ker\sigma_2\not\subset\ker\sigma_1$, otherwise we would have
		\[2 = [\Gamma:\ker\sigma_2] = [\Gamma:\ker\sigma_1]\cdot[\ker\sigma_1:\ker\sigma_2]= 2\cdot [\ker\sigma_1:\ker\sigma_2],\]
	which implies $\ker\sigma_1 = \ker\sigma_2$. Similarly, $\ker\sigma_1\not\subset\ker\sigma_2$. Thus there exists $\gamma\in\Gamma$ such that $\gamma\in\ker\sigma_2$ and $\gamma\notin\ker\sigma_1$, that is, $\gamma\in\ker\sigma_2\cap\left(\Gamma\backslash\ker\sigma_1\right)$. In this case $\gamma$ is AE. There also exists $\delta\in\Gamma$ such that $\delta\in\ker\sigma_1\cap\left(\Gamma\backslash\ker\sigma_2\right)$. Therefore $\delta$ is SR. Hence, $\sigma_1(\gamma\delta) = \sigma_1(\gamma)\sigma_1(\delta) = -1$ and $\sigma_2(\gamma\delta) = \sigma_2(\gamma)\sigma_2(\delta) = -1$, that is,  $\gamma\delta$ is AR. \cqd 

Consider the homomorphism $\sigma_1\sigma_2:\Gamma\to{\mathbb{Z}}_2$ defined by $(\sigma_1\sigma_2)(\gamma) = \sigma_1(\gamma)\sigma_2(\gamma)$, for all $\gamma\in\Gamma$. If $\sigma_1$ and $\sigma_2$ are non-trivial and distinct, then the product $\sigma_1\sigma_2$ is a third epimorphism, which is distinct from the both ones. 


We remark that the hypothesis $H(p) = 0$ is not necessary for the reciprocal of the next theorem.
\begin{theorem}\label{teogenebuzzi1}Let $\Gamma$ be a Lie group acting ${\sigma_1}$-semisymplectically on a symplectic vector space $(V,\omega)$ and $X_H$ an $\omega$-Hamiltonian vector field on $V$, with $H:V\to{\mathbb{R}}$ satisfying $H(p) = 0$ for some $p\in V$. If $X_H$ is $\Gamma_{\sigma_2}$-equivariant, then the Hamiltonian function $H$ is $\Gamma_{\sigma_1\sigma_2}$-invariant. Reciprocally, if $H$ is $\Gamma_{\sigma_2}$-invariant, then $X_H$ is $\Gamma_{\sigma_1\sigma_2}$-equivariant.
\end{theorem}       
\Dem Suppose $X_H(\rho_{\gamma}(x)) = \sigma_2(\gamma)d(\rho_{\gamma})_x X_H(x)$, for all $\gamma\in\Gamma$ and $x\in V$. Since $X_H$ is $\omega$-Hamiltonian, we have for all $x, \mathbf{v}\in V$ that 
	\begin{align}
		dH_x\mathbf{v}&=\omega(X_H(x),\mathbf{v}) =\sigma_1(\gamma)\omega(d(\rho_{\gamma})_x X_H(x),d(\rho_{\gamma})_x \mathbf{v})\nonumber\\
		& =\sigma_1(\gamma)\omega(\sigma_2(\gamma)X_H(\rho_{\gamma} (x)), d(\rho_{\gamma})_x\mathbf{v}) 	=(\sigma_1\sigma_2)(\gamma)dH_{\rho_{\gamma}(x)}d(\rho_{\gamma})_x\mathbf{v}\nonumber\\
		&=d((\sigma_1\sigma_2)(\gamma)(H\circ\rho_{\gamma}))_x\mathbf{v},\nonumber
	\end{align} where the second equality holds since the action of $\Gamma$ on $(V,\omega)$ is $\sigma_1$-semisymplectic. Thereby $(\sigma_1\sigma_2)(\gamma)(H\circ\rho_{\gamma})(x) = H(x)+c$, for some constant $c\in{\mathbb{R}}$, for all $\gamma\in\Gamma$ and $x\in V$. In particular, by considering $\gamma = 1$ the identity of $\Gamma$, we have $c = H(p)+c = (\sigma_1\sigma_2)(1) H(\rho_{1}(p)) = H(p) = 0.$ Therefore, $H(\rho_{\gamma}(x)) = (\sigma_1\sigma_2)(\gamma)H(x)$.
	
	Conversely, given $\gamma\in\Gamma$ and $x, \mathbf{v}\in V$, since $H(\rho_{\gamma}(x)) = \sigma_2(\gamma) H(x)$, we have $\sigma_2(\gamma)dH_x\mathbf{v} = d(H\circ {\rho}_{\gamma})_x\mathbf{v} = dH_{{\rho}_{\gamma}(x)}d({\rho}_{\gamma})_x\mathbf{v}$. Thus
	\begin{align}
		\omega(X_H({\rho}_{\gamma}(x)),d({\rho}_{\gamma})_x\mathbf{v})&=dH_{{\rho}_{\gamma}(x)} d({\rho}_{\gamma})_x\mathbf{v} = \sigma_2(\gamma)dH_x\mathbf{v}\nonumber\\
		& =\sigma_2(\gamma) \omega(X_H(x),\mathbf{v}) = \sigma_2(\gamma)\sigma_1(\gamma)\omega(d({\rho}_{\gamma})_x X_H(x),d({\rho}_{\gamma})_x \mathbf{v})\nonumber\\
		&=\omega((\sigma_1\sigma_2)(\gamma)d({\rho}_{\gamma})_x X_H(x),d({\rho}_{\gamma})_x\mathbf{v}),\nonumber
	\end{align} where the last equality follows by the bilinearity of $\omega$. Therefore,
 \[\omega(X_H(\rho_{\gamma}(x))-(\sigma_1\sigma_2)(\gamma) d(\rho_{\gamma})_xX_H(x), d(\rho_{\gamma})_x\mathbf{v}) = 0,\]
for all $\gamma\in\Gamma$ and $x, \mathbf{v}\in V$. Since $\omega$ is non-degenerate and $d(\rho_{\gamma})_x$ is an isomorphism, we conclude that $X_H(\rho_{\gamma}(x)) = (\sigma_1\sigma_2)(\gamma)d(\rho_{\gamma})_xX_H(x)$, for all $\gamma\in\Gamma$ and $x\in V$.\cqd

The last two results can be naturally extended to symplectic manifolds. As consequence of Theorem \ref{teogenebuzzi1} we present the following table, which also appears in \cite{Alomair} for the case where $\Gamma$ acts orthogonally on $({\mathbb{R}}^{2n},\omega_0)$. The signs of $\omega$, $X_H$ and $H$ represent, respectively, the signs of $\sigma_1(\gamma)$, $\sigma_2(\gamma)$ and $(\sigma_1\sigma_2)(\gamma)$ for $\gamma\in\Gamma$ satisfying the condition in the first column.

	\begin{table}[!htb] 
	\begin{center} 
		\begin{tabular}{|l|c|c|c|}
			\hline
			\hspace{2cm} Symmetry&$\omega$&$X_H$&$H$\\ \hline
			Symplectic equivariant (SE)&$+1$&$+1$&$+1$\\ \hline
			Symplectic reversible (SR)&$+1$&$-1$&$-1$\\ \hline
			Antisymplectic equivariant (AE)&$-1$&$+1$&$-1$\\ \hline
			Antisymplectic reversible (AR)&$-1$&$-1$&$+1$\\\hline
		\end{tabular} 
		\caption{\label{tabelatipos}Possible types of symmetry in the Hamiltonian context.}
		\end{center}
		\end{table}


In \cite{buono, Montaldi} the authors consider a unique homomorphism $\sigma: = \sigma_1=\sigma_2$ so that a symmetric Hamiltonian vector field admits only symmetries of type SE or AR. The next example presents a Hamiltonian vector field with the four types of symmetry of Table \ref{tabelatipos}.

\begin{example}\label{exemploum}{\rm
		Consider $\Gamma = {\mathbb{Z}}_2^{\tau}\times{\mathbb{Z}}_2^{\psi}$  acting on $(\mathbb{R}^4,\omega)$ by matrix multiplication, where\[\tau =\left[\begin{array}{cccc}
			-1&0&0&0\\
			0&1&0&0\\
			0&0&-1&0\\
			0&0&0&1
		\end{array}\right],\quad \psi =\left[\begin{array}{cccc}
			-1&0&0&0\\
			0&1&0&0\\
			0&0&1&0\\
			0&0&0&-1
		\end{array}\right]\quad\text{and}\quad [\omega] = \left[\begin{array}{cccc}
			0&0&0&-1\\
			0&0&1&0\\
			0&-1&0&0\\
			1&0&0&0
		\end{array}\right].\] Notice that $\psi\in Sp_{\omega}(2;{\mathbb{R}})$ and $\tau\in Sp_{\omega}^{-1}(2;{\mathbb{R}})$. Then $\Gamma$ is a subgroup of $\Omega_2$ and, by Theorem \ref{acao subgrupo de Omega}, its action on $(\mathbb{R}^4,\omega)$ is $\sigma_1$-semisymplectic so that $\ker\sigma_1 = \left\{(I_4,I_4),(I_4,\psi)\right\}$. Given $\lambda\in\mathbb{R}^*$, the $\omega$-Hamiltonian vector field \[X(x_1,x_2,x_3,x_4) = \lambda(x_2,-x_1,x_4,-x_3)\] associated with $H(x_1,x_2,x_3,x_4) = -\lambda\left(x_1x_3+x_2x_4\right)$ is $\Gamma_{\sigma_2}$-equivariant, with $\ker\sigma_2 = \left\{(I_4,I_4), (\tau, \psi)\right\}$. By Theorem \ref{teogenebuzzi1}, $H$ is $\Gamma_{\sigma_1\sigma_2}$-invariant, where $\sigma_1\sigma_2:\Gamma\to{\mathbb{Z}}_2$ is the epimorphism such that $\ker\sigma_1\sigma_2 = \left\{(I_4,I_4), (\tau, I_4)\right\}$. By Proposition \ref{proptipos}, the vector field $X$ presents four types of symmetry: 
	\begin{table}[htb]
			\centering
				\begin{tabular}{|c|c|c|c|}\hline
				SE&SR&AE&AR\\ \hline
				$(I_4, I_4)$&$(I_4, \psi)$&$(\tau, \psi)$&$(\tau, I_4)$\\ \hline
			\end{tabular}
			\caption{Types of symmetry of $X(x_1,x_2,x_3,x_4) = \lambda(x_2,-x_1,x_4,-x_3)$.}
					\label{tabela4tipos}
		\end{table}

}\end{example}

\section{Normal forms of $\omega$-Hamiltonian vector fields}\label{secaoformasnormaisNOVA}

The aim of this section is to present the normal form theory of $\omega$-Hamiltonian vector fields on $(V, \omega)$, based on the approach given in \cite{Chow} for the canonical symplectic space $(\Re^{2n}, \omega_0)$. Moreover, we establish a relation between the normal forms of an $\omega$-Hamiltonian vector field and of its associated Hamiltonian function. 

For what follows, we denote respectively by ${\mathcal P}^k$ and ${\overrightarrow{{\mathcal{P}}}^{k}}$ the vector spaces of the homogeneous maps of degree $k$ in ${\mathcal P}$ and in $\overrightarrow{{\mathcal{P}}}$. Similarly, ${\mathcal P}^k(\Gamma)$ and ${\overrightarrow{{\mathcal{P}}}^{k}}(\Gamma)$ are the vector spaces of the homogeneous maps of degree $k$ in ${\mathcal P}(\Gamma)$ and in $\overrightarrow{{\mathcal{P}}}(\Gamma)$ (see (\ref{conjSIMETRIAS})), respectively.

Consider $(V,\omega)$ a symplectic vector space of dimension $2n$ and $X (x) = ([\omega]^{-1})^T\nabla H(x)$ an $\omega$-Hamiltonian vector field on $V$ with an isolated singularity at the origin, that is, $X(0) = 0$. Hence we have $\nabla H(0) = 0$ and assume $H(0) = 0$. Thus, the Taylor series of $H$ centered at the origin is written as\begin{equation}
	\label{serieH}
	H^2(x)+H^3(x)+\cdots+H^r(x)+O(|x|^{r+1}),
\end{equation}where $H^k\in{\mathcal{P}}^k$. 

\begin{lemma}\label{LH2}The $\omega$-Hamiltonian vector field $X_{H^2}$ associated with $H^2$ determines the linear part of $X$, that is, $X_{H^2}(x) = Lx,$ where $L = dX_0$ is the linearization of $X$ at the origin.
\end{lemma}
\Dem Since $X(x) = ([\omega]^{-1})^T\nabla H(x) = (dH_x[\omega]^{-1})^T$, we obtain\[L = dX_0 = (d^2H_0[\omega]^{-1})^T = ([\omega]^{-1})^Td^2H_0.\] By definition, we have $H^2(x) = \frac{1}{2}\mathbf{x}^Td^2H_0\mathbf{x}$. Thus $\nabla H^2(x) = d^2H_0x$ and\[X_{H^2}(x) = ([\omega]^{-1})^T\nabla H^2(x) = ([\omega]^{-1})^Td^2H_0x = Lx.\] \cqd 

Based on the classical normal form theory, we want to determine a change of coordinates $x = \xi (y)$ that preserves the linearization $L = dX_0$ and, naturally, the dynamics of the system $\dot{x} = X(x)$ around the origin. For each $k\geq3$,\[X_{H^2+\cdots+H^k}(x) = X_{H^2}(x)+\cdots+X_{H^k}(x) = Lx+X_{H^3}(x)+\cdots+X_{H^k}(x),\]where $X_{H^k} = ([\omega]^{-1})^T\nabla H^k\in{{\overrightarrow{{\mathcal{P}}}^{k-1}}}$. From now on, we assume that\begin{equation}\label{1.1 linha}X(x) \ = Lx+X_{H^3}(x)+X_{H^4}(x)+\cdots+X_{H^r}(x)+O(|x|^r)
\end{equation}and denote $\mu(x) = X_{H^3}(x)+X_{H^4}(x)+\cdots+X_{H^r}(x)+O(|x|^{r})$. Consider $x=\xi(y)$ a change of coordinates in a neighborhood $U$ of the origin, that is, $d\xi_y$ is invertible with $\xi(0)=0.$ Thus \begin{equation}\label{1.2}
	\xi^{-1}_*X(y) \ = \ d\xi_y^{-1}(L\xi(y)+\mu(\xi(y))), \end{equation} with $y \in U$. We denote by $\varphi(y)$ the right side of (\ref{1.2}) and by $\nu$ the derivative of $d\xi_y^{-1}$ with respect to $y$. Then the Jacobian of $\varphi$ is given by\[J_{\varphi}(y)=\nu(L\xi(y)+\mu(\xi(y)))+d\xi_y^{-1}(Ld\xi_y+d\mu_{\xi(y)}d\xi_y).\]For $y=0$, we have $J_{\varphi}(0) = d\xi_0^{-1}Ld\xi_0+d\xi_0^{-1}d\mu_0d\xi_0 = d\xi_0^{-1}Ld\xi_0,$ since $d\mu_0=0$. Thus, the linear part of (\ref{1.2}) is $d\xi_0^{-1}Ld\xi_0y$.

Assume that $\xi$ has the form $\xi(y)=y+O(|y|^2),$ with $y\rightarrow0$. Then $d\xi_0 = d\xi_0^{-1} = {\rm Id}$ whence the linear part of (\ref{1.2}) is $Ly$. Therefore, \begin{equation}\label{1.3}
	\xi^{-1}_*X(y) \ = Ly+g(y),
\end{equation}with $y\in U$, where $g(y)$ has only terms of degree greater than or equal to 2.

The idea is determining $\xi$ so that (\ref{1.3}) is still an $\omega$-Hamiltonian vector field. This is possible if  $\xi$ is a symplectic diffeomorphism and, in this case, we obtain a new $\omega$-Hamiltonian vector field associated with $H \circ \xi$ (see \cite[Theorem 4.2]{BRHT2021}). 


\begin{proposition}\label{mudancasimpletica}For each $k\geq3$, given $\xi^k\in{\mathcal{P}}^k$, there exist a neighborhood of the origin $U_k\subset V$  and a symplectic diffeomorphism $\xi:U_k\to V$ written as\begin{equation}
		\label{mudsimpletica}
		\xi(y) = y + X_{\xi^k}(y)+ O(|y|^{k}),
	\end{equation}with $y\in U_{k}$, where $X_{\xi^k} = ([\omega]^{-1})^T\nabla\xi^k$.
\end{proposition}
\Dem Consider the dynamic system $\phi$ associated with $Y = X_{\xi^k}$. Hence $\phi(t,y)$ satisfies the initial value problem\begin{equation}\label{pvi}
		\begin{cases}
			\dot{\phi}(t,y)= Y(y) \\   \phi(0,y) = y   
		\end{cases}
	\end{equation}and it can be expanded in a Taylor series around the origin as\begin{equation}
		\label{expansaosimpletica}
		\phi(t,y) = \phi^1(t)y+\phi^2(t,y)+\cdots+\phi^{k-1}(t,y)+O(|y|^k),
	\end{equation}where $\phi^1:\mathbb{R}\to{\mathbb{M}}_{2n}(\mathbb{R})$, $\phi^j:\mathbb{R}\times V\to V$ and, for all fixed $t$, $\phi^j(t,\cdot)\in{{\overrightarrow{\mathcal{P}}^{j}}}$, $j=2,\ldots,k-1$.
	Since $Y$ is a homogeneous polynomial map of degree $k-1$, by replacing (\ref{expansaosimpletica}) in (\ref{pvi}) we get
	\begin{center}
		$\begin{cases}
			\dot{\phi}^1(t) = 0 \\
			\phi^1(0)= I_ {2n}
		\end{cases}$\quad  and\qquad $\begin{cases}
			\dot{\phi}^j(t,y) = 0 \\
			\phi^j(0,y)= 0
		\end{cases}$,
	\end{center}for all $j=2,\ldots,k-2$.  Thus, given $2\leq j\leq k-2$, for all $t\in\mathbb{R}$ in the maximal interval and all $y\in V$ we have $\phi^1(t) = I_{2n}$ and $\phi^j(t,y) = 0$, since the first line in each IVP ensures us that $\phi^1$ and $\phi^j(\cdot,y)$ are constant in $t$. Moreover\begin{equation}\nonumber
		\begin{cases}
			\dot{\phi}^{k-1}(t,y) = Y(y) \\
			\phi^{k-1}(0,y)=0
		\end{cases}.
	\end{equation}Now, by integrating with respect to $t$ and using the initial condition we have $\phi^{k-1}(t,y) = tY(y)$. Hence $\phi^{k-1}(1,y) = Y(y) = X_{\xi^k}(y)$ and the map $Y_1$ defined by $Y_1(y) = \phi(1,y)$ is symplectic (see \cite[Proposition V.24]{Zehnder}). By (\ref{expansaosimpletica}),\[Y_1(y) = y+X_{\xi^k}(y)+O(|y|^k)\]is a symplectic map. Finally take  $\xi = Y_1$. By inverse mapping theorem, there exists a neighborhood $U_k\subset V$ such that $\xi$ is a local diffeomorphism. \cqd 

We now recall some concepts and results about Hamiltonian vector fields. Given two $C^{\infty}$ functions $F, G: V \to \Re$, the Poisson bracket of $F$ and $G$ is a function $\{F,G\}: V \to \Re$ defined by 
$$\{F,G\}(x) = \omega \left(X_F(x),X_G(x)\right),$$ where $X_F$ and $X_G$ are $\omega$-Hamiltonian vector fields associated with $F$ and $G$, respectively. It is possible to prove that $\{F,G\}(x) = \langle\nabla F(x), X_G(x)\rangle_{\rm can}$ for all $x\in V$, where $\langle\cdot,\cdot\rangle_{\rm can}$ is the canonical inner product on $V$. Moreover
\begin{equation}
\label{eq:poisson}
X_{\{F,G\}} = [X_G,X_F],  
\end{equation} where $[ \ , \ ]$ denotes the Lie bracket of vector fields. For more details, see \cite[Proposition V.28]{Zehnder}).


Fixed $k\geq3$, consider a symplectic change of coordinates near the identity as in (\ref{mudsimpletica}). For each $y\in U_{k}$, we have $d\xi_y={\rm Id} + d(X_{\xi^{k}})_y+O(|y|^{k-1})$, whence $d\xi_y^{-1} ={\rm Id}-d(X_{\xi^{k}})_y+O(|y|^{k-1})$. So we can rewrite (\ref{1.2}) like

\begin{eqnarray}\xi^{-1}_*X(y)&=&d\xi_y^{-1}(L\xi(y)+\mu(\xi(y)))\nonumber\vspace{0.3cm}\\
	&=&\left({\rm Id}-d(X_{\xi^{k}})_y+O(|y|^{k-1})\right)\left(Ly+LX_{\xi^{k}}(y)+ \mu(\xi(y)) \right)\nonumber\vspace{0.3cm} \\
	&=&Ly+LX_{\xi^{k}}(y)+X_{H^3}(y)+\cdots+X_{H^{k}}(y)-d(X_{\xi^{k}})_yLy\nonumber\vspace{0.1cm}\\&&-d(X_{\xi^{k}})_yLX_{\xi^{k}}(y)-d(X_{\xi^{k}})_y(\mu(\xi(y))) + O(|y|^k)\nonumber\vspace{0.3cm}\\
	&=&Ly+X_{H^3}(y)+\cdots+X_{H^{k-1}}(y)\nonumber\vspace{0.1cm}\\&&+\left(X_{H^{k}}(y)-(d(X_{\xi^{k}})_yLy-LX_{\xi^{k}}(y))\right)+O(|y|^{k}),\nonumber
\end{eqnarray}with $y\in U_{k}$. Since $d(X_{\xi^{k}})_yL-LX_{\xi^{k}}$ is the Lie bracket $[L, X_{\xi^k}]$, it follows by Lemma \ref{LH2} and by (\ref{eq:poisson}) that the term in  parentheses is equal to\[X_{H^{k}}(y)-[X_{H^2},X_{\xi^{k}}](y) = X_{H^{k}}(y)-X_{\{\xi^{k},H^2\}}(y) = X_{H^{k}-\{\xi^{k},H^2\}}(y).\] Therefore, the vector field (\ref{1.1 linha}) is formally conjugated\footnote{A vector field $Y$ is said {formally conjugated} to $X$ if a power series of $Y$ is conjugated to a power series of $X$, as formal vector fields.} to
\begin{equation}\label{1.5}
	\xi^{-1}_*X(y) = Ly+X_{H^3}(y)+\cdots+X_{H^{k-1}}(y)+X_{H^{k}-\{\xi^{k},H^2\}}(y)+O(|y|^{k}),
\end{equation} in which the linear part $L$ and the terms of degree less than $k-1$ of (\ref{1.1 linha}) are preserved.


\begin{definition}Given a quadratic form $H^2:V\to\mathbb{R}$, for each $k\geq3$ we define the linear operator $ad_{H^2}^k: \ {\mathcal{P}}^k\to{\mathcal{P}}^k$ by \begin{equation}\label{homologicokh2} ad_{H^2}^k(\xi^{k})= \{\xi^k,H^2\}.
	\end{equation}
\end{definition}

For the next theorem, consider ${\mathcal{D}}^{k}$ a complementary subspace of $ad_{H^2}^k({\mathcal{P}}^k)$ in ${\mathcal{P}}^k$, that is,
\begin{equation}\label{1.6}
	{\mathcal{P}}^k = ad_{H^2}^k({\mathcal{P}}^k)\oplus{\mathcal{D}}^{k}.
\end{equation}

\begin{theorem}\label{teo1.1} Let $(V,\omega)$ be a symplectic vector space and $X:V\to V$ an $\omega$-Hamiltonian vector field such that $X(0) = 0$ and $dX_0=L$. Suppose that the Hamiltonian function associated with $X$ admits a Taylor expansion at the origin as in (\ref{serieH}) and consider the decomposition (\ref{1.6}), for $k=3,\ldots,r$. Then $X$ is formally conjugated to the $\omega$-Hamiltonian vector field
	\begin{equation}\label{1.14}
		Y(y) = Ly+X_{K^3}(y)+\cdots+X_{K^r}(y)+O(|y|^{r})
	\end{equation} associated with $K = H^2+K^3+\cdots+K^r + O(|y|^{r})$, where $K^k\in{\mathcal{D}}^k$ for $k=3,\ldots,r$.
\end{theorem}
\Dem The idea is to exhibit a sequence of symplectic changes of coordinates near the identity of the form $x=y+X_{\xi^{k}}(y)+O(|y|^{k}),$ with $y\in U_{k}$, where ${\xi^{k}}\in{\mathcal{P}}^{k}$ and $U_{k}$ is a neighborhood of the origin, with $U_{k+1}\subseteq U_k$, such that in the new coordinates the vector field (\ref{1.1 linha}) is written as in (\ref{1.14}).
	
	By (\ref{1.6}), there exist $G^3\in ad_{H^2}^3({\mathcal{P}}^3)$ and $K^3\in{\mathcal{D}}^3$ such that ${H^3}=G^3+K^3$. Hence we can choose $\xi^3\in{\mathcal{P}}^3$ so that $G^3= ad_{H^2}^3(\xi^3)=\{\xi^3,H^2\}$. Consider the symplectic diffeomorphism given by $\xi(y) = y+ X_{\xi^3}(y)+O(|y|^3),$ with $y\in U_3$, for some neighborhood $U_3$ of the origin. Then $X$ is conjugated to $\xi^{-1}_*X$ as in (\ref{1.5}), that is, $$\xi^{-1}_*X(y) =Ly+X_{H^3-ad_{H^2}^3(\xi^{3})}(y)+O(|y|^3) = Ly+X_{K^3}(y)+O(|y|^3).$$ Now we proceed by induction on $k$, by remembering that such change of coordinates do not modify terms of degree less than $k-1$. We assume for $3\leq k -1<r$ the existence of a sequence of symplectic diffeomorphisms as in (\ref{mudsimpletica}) such that $X$ is conjugated to
	\begin{equation}\label{1.5 na prova}
		\xi^{-1}_*X(y) = Ly+X_{K^3}(y)+\cdots+X_{K^{k-1}}(y)+O(|y|^{k-1}),
	\end{equation}with $y\in U_{k-1}$, where $U_{k-1}\subset U_{k-2}$ is a neighborhood of the origin and $K^j\in{\mathcal{D}}^j$ for $j=3,\ldots,k-1$. Since (\ref{1.5 na prova}) is an $\omega$-Hamiltonian vector field we have\begin{equation}\label{1.5novo}                                                                                                                                                                                                                                                            
		\xi^{-1}_*X(y) = Ly+X_{K^3}(y)+\cdots+X_{K^{k-1}}(y)+X_{{\overline{H}}^{k}}(y)+X_{{\overline{H}}^{k+1}}(y)+\cdots,
	\end{equation}where $X_{{\overline{H}}^{k+j}}$ is a homogeneous polynomial map of degree $k+j-1$ for $j\geq0$.
	
	Consider $\xi(y) =  y+X_{\xi^{k}}(y)+O(|y|^{k})$, with $y\in U_{k}$, for some neighborhood of the origin $U_{k}\subseteq U_{k-1}$, where we choose $\xi^{k}\in{\mathcal{P}}^{k}$ such that ${\overline{H}}^{k}=ad_{H^2}^{k}(\xi^{k})+K^{k}$ for some $K^{k}\in{\mathcal{D}}^{k}$. Write $Y(y): = \xi^{-1}_*X(y)$. Hence (\ref{1.5novo}) is given by \begin{align*}
		Y(y)&=Ly+X_{K^3}(y)+\cdots+X_{K^{k-1}}(y)+X_{{\overline{H}}^{k}-ad_{H^2}^k(\xi^k)}(y)+O(|y|^k)\\
		&=Ly+X_{K^3}(y)+\cdots+X_{K^{k-1}}(y)+X_{K^k}(y)+O(|y|^k),
	\end{align*} with $y\in U_k$. \cqd
	

The previous result induces the following definition.
\begin{definition}\label{fnfuncaodefi}Given an Hamiltonian function $H:V\to\mathbb{R}$ such that $H(0)=0$ and $\nabla H(0) = {0}$, a normal form of $H$ up to order $r$ is a sum
	$$		K = H^2+K^3+\cdots+K^r,$$ where $H^2$ is the quadratic form associated with the Hessian of $H$ and $K^k\in{{\mathcal{D}}^k}$ for $k=3,\ldots,r$. If $X$ is an $\omega$-Hamiltonian vector field associated with $K$, we say that $$Y(y) = Ly+X_{K^3}(y)+\cdots+X_{K^r}(y)$$ is a normal form of $X$ up to order $r-1$.
\end{definition}

\quad 

In \cite{Be,Be2}, Belitskii develops a theory of normal forms for arbitrary vector fields $X: \mathbb{R}^n \to \mathbb{R}^n$ with $X(0) = 0$ and $L = dX_0$. For each $k\geq3$, the author defines the linear operator $Ad_L^{k-1}: {\overrightarrow{{\mathcal{P}}}^{k-1}}\to{\overrightarrow{{\mathcal{P}}}^{k-1}}$ as $Ad_L^{k-1}(\xi) = [L,\xi]$ and consider a subspace ${\mathcal{C}}^{k-1}$ in the decomposition $${{\overrightarrow{{\mathcal{P}}}^{k-1}}} = Ad_L^{k-1}({{\overrightarrow{{\mathcal{P}}}^{k-1}}})\oplus{\mathcal{C}}^{k-1}.$$ It is possible to prove that $X$ is formally conjugated to $Y(y) \ = Ly+g^2(y)+\cdots+g^r(y)+O(|y|^{r+1}),$ where $g^{k-1}\in{\mathcal{C}}^{k-1}$ (see \cite[Chapter 2, Theorem 1.1]{Chow}). From \cite[Chapter 2, Theorem 1.7]{Chow} we can assume ${\mathcal{C}}^{k-1} = \ker Ad_{L^T}^{k-1}$ and by \cite[XVI Lemma 5.4]{Golu2} we have 
\begin{equation}
\label{propmetodoElphick}
	\ker Ad_{L^T}^{k-1} = {{\overrightarrow{{\mathcal{P}}}^{k-1}}}(\mathbf{S}),
\end{equation} where $\mathbf{S}$ is defined in (\ref{grupoS}). In the next subsection, we will relate the normal form of an $\omega$-Hamiltonian vector field $X$ given in (\ref{1.14}) with the obtained by the classical theory developed in \cite{Be,Be2} for an appropriate subspace $\mathcal{D}^k$. 

\subsection{An algebraic method to determine ${\mathcal{D}}^k$}\label{secaoDK}

In this subsection, we exhibit a convenient complementary subspace ${\mathcal{D}}^k$ as defined in (\ref{1.6}). In \cite[Theorem 2.1, Lemma 2.3]{Montaldi}, the authors consider the canonical symplectic space $({\mathbb{R}}^{2n},\omega_0)$ and show that a possible choice for ${\mathcal{D}}^k$ is the set ${{\mathcal{P}}^k}({\bf S})$ of all \textbf{S}-invariant homogeneous polynomial functions of degree $k$. For this result, they use that ${\bf S}$ acts symplectically on $(\Re^{2n}, \omega_0)$. However, this property is not always true on an arbitrary symplectic vector space $(V,\omega)$ (see Example \ref{Snaosimpleticoex} and Theorem \ref{acao subgrupo de Omega}). 


From now on, we assume that the action of $\mathbf{S}$ on $(V, \omega)$ is symplectic. In order to characterize the subspace ${\mathcal{D}}^k$, we introduce an inner product in ${\mathcal{P}}^k$: given $F, G\in{\mathcal{P}}^k$ and $X_F$, $X_G$ their respective $\omega$-Hamiltonian vector fields, we define \[\ll F,G\gg = \langle X_F, X_G\rangle,\] where $\langle\cdot,\cdot\rangle$ is a convenient inner product in ${{\overrightarrow{{\mathcal{P}}}^{k-1}}}$ such that the equality (\ref{propmetodoElphick}) holds (see \cite[XVI Theorem 5.3]{Golu2}). Notice that a possible choice for ${\mathcal{D}}^k$ is the orthogonal complementary subspace of $ad_{H^2}^k({{\mathcal{P}}^k})$ characterized by $\ker (ad_{H^2}^k)^*$, where $(ad_{H^2}^k)^*$ is the adjoint operator of $ad_{H^2}^k$. Since ${\bf S}$ acts symplectically on $(V, \omega)$, we highlight that $Y(x) = L^Tx$ is an $\omega$-Hamiltonian vector field  (see Proposition \ref{exacaosimpleticasoS2} and \cite[Proposition 4.6]{BRHT2021}), that is, there exists a quadratic form denoted by $H^2_T: V\to\mathbb{R}$ such that $Y(x) = X_{H^2_{T}}(x)$. For each $k\geq3$, it is possible to prove that $ad_{H^2_{T}}^k$ is the adjoint operator of $ad_{H^2}^k$ with respect to the inner product $\ll \cdot,\cdot\gg$. The proof for this result is similar to that found in {\cite[Lemma 2.3]{Montaldi}} for $(\mathbb{R}^{2n},\omega_0)$. 

Therefore, in order to obtain a normal form up to order $r$ of a Hamiltonian function $H:V\to\mathbb{R}$, it is enough to determine $\ker ad_{H^2_{T}}^k$, $k=3,\ldots,r$. By properties of Poisson bracket and by Lemma \ref{LH2}, \begin{equation}\nonumber
	ad_{H^2}^k(\xi^{k})(x) = \{\xi^k,H^2\}(x) = \langle\nabla\xi^k(x), X_{H^2}(x)\rangle_{\rm can} = \langle\nabla\xi^k(x),Lx\rangle_{\rm can}.
\end{equation} Hence we need to solve the partial differential equation $\langle\nabla\xi^k(x),L^Tx\rangle_{\rm can} = 0,$ with $x\in V$. However we present an alternative algebraic characterization for this kernel. First, we remark that $F\in\ker ad^k_{H^2_T}$ if and only if $\{F,H^2_T\} = ad_{H^2_T}^k(F) = 0$, which occurs if and only if $X_{\{F,H^2_T\}} = ([\omega]^{-1})^T\nabla\{F,H^2_T\} = 0$. Hence

\begin{equation} \label{esta em Ck}
\begin{array}{rcl}
\ker ad_{H^2_{T}}^k & = & \{F\in{{\mathcal{P}}^k}: X_{\{F,H^2_{T}\}} = 0\} = \{F\in{{\mathcal{P}}^k}: [X_{H^2_{T}},X_F] = 0\} \vspace{.1cm} \\
	&= &\{F\in{{\mathcal{P}}^k}: X_F \in{\overrightarrow{{\mathcal{P}}}^{k-1}}({\bf S})\} = {\mathcal{P}}^{k}({\bf S}),
\end{array}
\end{equation} where the third equality follows from (\ref{propmetodoElphick}), since $Ad_{L^T}^{k-1}(X_F)= [X_{H^2_{T}},X_F]$. In the last equality we use that the action of $\textbf{S}$ on $(V,\omega)$ is symplectic and  Theorem \ref{teogenebuzzi1}, by considering $\Gamma = \textbf{S}$ and $\sigma_1=\sigma_2 :\mathbf{S}\to\mathbb{Z}_2$ the trivial homomorphism. 



Therefore, a good choice for the complementary subspace ${\mathcal{D}}^{k}$ in (\ref{1.6}) is the vector space ${\mathcal{P}}^{k}({\bf S})$, that is, 
\begin{equation}
	\label{complemento elphick}
	{\mathcal{P}}^{k} = ad^k_{H^2}({\mathcal{P}}^{k})\oplus {\mathcal{P}}^{k}({\bf S}).
\end{equation}Thus, in order to describe a normal form up to order $r-1$ of the $\omega$-Hamiltonian vector field $X_H$, it is enough to characterize the ring ${\mathcal{P}}^{k}({\bf S})$, for $k=3,\ldots,r$. In fact, we show in Theorem \ref{teo1.1} that the original vector field is formally conjugated to (\ref{1.14}). By (\ref{esta em Ck}), if $\mathcal{D}^k = {\mathcal{P}}^{k}({\bf S})$, then $X_{K^k}\in{\overrightarrow{{\mathcal{P}}}^{k-1}}({\bf S})$, where ${\overrightarrow{{\mathcal{P}}}^{k-1}}({\bf S})$ is a choice for the complementary subspace $\mathcal{C}^{k-1}$ according to (\ref{propmetodoElphick}). 

%

Therefore, if $L$ is a non-null matrix, by Theorem \ref{teo1.1} and by (\ref{complemento elphick}), the normal form of the $\omega$-Hamiltonian vector field $X$ presents symmetries, even if $X$ does not have any symmetry. Similarly for the normal form of the Hamiltonian function associated with $X$.

\section{Normal forms of $\omega$-Hamiltonian vector fields with symmetries}\label{secaoformasnormaissimetrias}
Fixed a compact Lie group $\Gamma$ acting linearly on a finite-dimensional real vector space, Golubitski, Stewart and Schaeffer \cite[XVI Theorem 5.8; Theorem 5.9]{Golu2} present a normal form of a $\Gamma$-equivariant vector field, based on method of Elphick et al. \cite{elphic}. Afterwards, Baptistelli, Manoel and Zeli {\cite[Theorem 4.7]{Bap}} generalize such a result in the context of $\Gamma_{\sigma}$-equivariant vector fields. The union of the classical method of normal forms in the Hamiltonian context presented in Section \ref{secaoformasnormaisNOVA} with the algebraic method developed in \cite{Bap} is the starting point for the study of normal forms of $\omega$-Hamiltonian vector fields with symmetries taking into account the semisymplectic action of $\Gamma$ on $(V,\omega)$.


From now on, $\Gamma$ is a compact Lie group and $\sigma_1,$ $\sigma_2:\Gamma\to{\mathbb{Z}}_2$ are two group homomorphisms. We assume that the action of $\Gamma$ on $(V,\omega)$ is linear and ${\sigma_1}$-semisymplectic, and $X:V\to V$ is a $\Gamma_{\sigma_2}$-equivariant $\omega$-Hamiltonian vector field, with $X(0) = 0$. We consider $H:(V,\omega)\to\mathbb{R}$ the Hamiltonian function associated with $X$ such that $H(0)=0$ and $\nabla H(0) = 0$. By Theorem \ref{teogenebuzzi1}, $H$ is $\Gamma_{\sigma_1\sigma_2}$-invariant. Our goal is presenting a normal form of $X$ which is also $\Gamma_{\sigma_2}$-equivariant and $\omega$-Hamiltonian. 


Let $L = dX_0$ be the linearization of $X$ at the origin and $\mathbf{S}$ the linear Lie group defined in (\ref{grupoS}). If $\sigma_2$ is an epimorphism, we have a well defined action of ${\bf S}\rtimes\Gamma$ on $V$, namely, \begin{equation}\label{doisasterisco}
	(e^{sL^T},\gamma)x = e^{sL^T}(\gamma x) = \gamma e^{\sigma_2(\gamma)sL^T}x,
\end{equation} as in \cite{Bap}. If $\sigma_2$ is trivial, then the actions of $\Gamma$ and $\mathbf{S}$  commute, since there exists an inner product on 
$V$ such that the action of $\Gamma$ is orthogonal (see \cite[XII Proposition 1.3]{Golu2}). In this case, we also assume the action in (\ref{doisasterisco}) by taking $\sigma_2(\gamma) = 1$ for all $\gamma \in \Gamma$. 

%


For what follows, we denote $\mathcal{P}_{\sigma}^k(\Gamma): = \mathcal{P}^k \cap \mathcal{P}_{\sigma}(\Gamma)$ (see (\ref{conjSIMETRIAS})). In order to obtain a normal form as in Theorem \ref{teo1.1}, we present now a result that provides us a complementary subspace of $ad_{H^2}^k(\mathcal{P}_{\sigma_1}^k(\Gamma))$ in $\mathcal{P}_{\sigma_1\sigma_2}^k(\Gamma)$, where $H^2$ is the quadratic form of $H$ and $ad_{H^2}^k$ is defined in (\ref{homologicokh2}). A proof of this result is presented in Section \ref{adsubsecaopropriedades}. 

\begin{theorem}\label{teo1.1nosso}
	If $\mathbf{S}$ acts symplectically on $(V,\omega)$, we have the decomposition
	\begin{equation}
		\label{soma nossa}
		{\mathcal{P}}^{k}_{\sigma_1\sigma_2}(\Gamma) = ad_{H^2}^{k}({\mathcal{P}}^{k}_{\sigma_1}(\Gamma))\oplus{\mathcal{P}}^{k}_{\tilde{\sigma}_1\tilde{\sigma}_2}\left({\bf S}\rtimes\Gamma\right),
	\end{equation}
	for each $k\geq3$, where $\tilde{\sigma}_j:{\bf S}\rtimes\Gamma\to{\mathbb{Z}}_2$ is defined by\begin{equation}\label{estrelaCorrecPatMuitasvezes}
		\tilde{\sigma}_j(e^{sL^T},\gamma) = \sigma_j(\gamma),
	\end{equation}with $j=1,2$, for all $s\in\mathbb{R}$ and $\gamma\in\Gamma$.
\end{theorem}

Our aim is to preserve the $\omega$-Hamiltonian condition and the symmetries of the original vector field. For this, we consider changes of coordinates $\xi$ which are symplectic, as those described in (\ref{mudsimpletica}), and $\Gamma$-equivariant, that is, $\xi\circ \rho_{\gamma} = \rho_{\gamma}\circ \xi,$ for all $\gamma\in\Gamma$. In this case, the pushforward $\xi^{-1}_*X$ is $\Gamma_{\sigma_2}$-equivariant, since \begin{equation}\nonumber
	(\rho_{\gamma})_{*}(\xi^{-1}_*X) =  (\rho_{\gamma}\circ\xi^{-1})_*X =(\xi^{-1}\circ\rho_{\gamma})_*X = \xi^{-1}_*((\rho_{\gamma})_{*}X) =  \sigma_2(\gamma)\xi^{-1}_*X,
\end{equation}for all $\gamma\in\Gamma$. 

In the next lemma, we present a condition for a symplectic diffeomorphism to be $\Gamma$-equivariant.


\begin{lemma}\label{mudsimpletica equivariante}Suppose $\Gamma$ acts $\sigma_1$-semisymplectically on $(V,\omega)$ and let $\xi$ be a symplectic diffeomorphism as in (\ref{mudsimpletica}). If $\xi^k\in{\mathcal{P}}^{k}_{\sigma_1}(\Gamma)$, then $\xi$ is $\Gamma$-equivariant.
\end{lemma}
\Dem By Theorem \ref{teogenebuzzi1}, if the action of $\Gamma$ on $(V,\omega)$ is $\sigma_1$-semisymplectic, then ${\xi^k}$ is $\Gamma_{\sigma_1}$-invariant if and only if $X_{\xi^k}$ is $\Gamma$-equivariant. Thus, $(\rho_{\gamma})_{*}X_{\xi^k} = X_{\xi^k}$ and since $\rho_{\gamma}:V\to V$ is a diffeomorphism we have\begin{equation}\nonumber
		\rho_{\gamma}\circ(X_{\xi^k})_t  =  ((\rho_{\gamma})_{*}X_{\xi^k})_t\circ\rho_{\gamma} = (X_{\xi^k})_t\circ\rho_{\gamma},
	\end{equation}for all $\gamma\in\Gamma$ and $t\in\mathbb{R}$, where $(X_{\xi^k})_t$ is the flow of $X_{\xi^k}$. In particular, for $t=1$ we have\begin{equation}\label{estrela6pontas}
		\rho_{\gamma}\circ(X_{\xi^k})_1 = (X_{\xi^k})_1\circ\rho_{\gamma}.
	\end{equation}In the proof of Proposition \ref{mudancasimpletica}, we choose $\xi$ as a solution for $\dot{y} = X_{\xi^k}(y)$. More precisely we take $\xi = (X_{\xi^k})_1$. By (\ref{estrela6pontas}) we conclude that $\rho_{\gamma}\circ \xi=\xi\circ \rho_{\gamma}$ for all $\gamma\in\Gamma$, that is, $\xi$ is $\Gamma$-equivariant. \cqd 


In the next theorem we present a normal form of an $\omega$-Hamiltonian and $\Gamma_{\sigma_2}$-equivariant vector field $X$ which preserves the $\omega$-Hamiltonian condition and the symmetries of $X$.

\begin{theorem}\label{fnormalnossa}Let $\Gamma$ be a compact Lie group and $\sigma_1,$ $\sigma_2:\Gamma\to\mathbb{Z}_2$ group homomorphisms. Suppose that the action of $\Gamma$ on $(V,\omega)$ is $\sigma_1$-semisymplectic and consider $X:V\to V$ an $\omega$-Hamiltonian and $\Gamma_{\sigma_2}$-equivariant vector field such that $X(0) = 0$ and $L = dX_0$. If $\mathbf{S}$ acts symplectically on $(V,\omega)$, then $X$ is formally conjugated to\begin{equation}\nonumber
		Y(y) = Ly+X_{K^3}(y)+\cdots+X_{K^{r}}(y)+O(|y|^{r}),
	\end{equation}where $K^{k}\in {\mathcal{P}}^{k}_{\tilde{\sigma}_1\tilde{\sigma}_2}\left({\bf S}\rtimes\Gamma\right)$, $k=3,\ldots,r$, and $\tilde{\sigma}_j:{\bf S}\rtimes\Gamma\to{\mathbb{Z}}_2$ is defined in (\ref{estrelaCorrecPatMuitasvezes}).
\end{theorem}

\Dem Consider $X$ as in (\ref{1.1 linha}). By $\Gamma_{\sigma_2}$-equivariance of $X$, we have that $L$ and $X_{H^k}$ are also $\Gamma_{\sigma_2}$-equivariant for all $k\geq3$. Thus, $X_{H^k}\in{\overrightarrow{\mathcal{P}}_{\sigma_2}^k(\Gamma)}$ and by Theorem \ref{teogenebuzzi1} we know that $H^k\in{\mathcal{P}_{\sigma_1\sigma_2}^k}(\Gamma)$ for each $k\geq3$. From here, the proof of this theorem follows similarly to the proof of the Theorem \ref{teo1.1}. In the induction hypothesis, we assume that for $3\leq k -1<r$ there exists a sequence of symplectic changes of coordinates which transforms (\ref{1.1 linha}) to the $\omega$-Hamiltonian and $\Gamma_{\sigma_2}$-equivariant vector field
	\begin{equation}\label{1.5novoDois}                                                                	Y(y) = Ly+X_{K^3}(y)+\cdots+X_{K^{k-1}}(y)+X_{{\overline{H}}^{k}}(y)+X_{{\overline{H}}^{k+1}}(y)+\cdots,
	\end{equation}where $K^j\in{\mathcal{P}}^{j}_{\tilde{\sigma}_1\tilde{\sigma}_2}\left({\bf S}\rtimes\Gamma\right)$ for $j=3,\ldots,k-1$. For each $k\geq3$ we have $X_{{\overline{H}}^{k}}\in{\overrightarrow{\mathcal{P}}_{\sigma_2}^k(\Gamma)}$ and, again by Theorem \ref{teogenebuzzi1}, ${\overline{H}}^k\in{\mathcal{P}}^{k}_{\sigma_1\sigma_2}(\Gamma)$. Hence, by (\ref{soma nossa}), there exist $G^k\in ad_{H^2}^k(\mathcal{P}^k_{\sigma_1}(\Gamma))$ and $K^k\in{\mathcal{P}}^{k}_{\tilde{\sigma}_1\tilde{\sigma}_2}\left({\bf S}\rtimes\Gamma\right)$ such that ${{\overline{H}}^k}=G^k+K^k$. Let $\xi^k\in{\mathcal{P}}^k_{\sigma_1}(\Gamma)$ such that $ad_{H^2}^k(\xi^k)=G^k$ and consider the symplectic diffeomorphism $\xi(y) =  y+X_{\xi^{k}}(y)+O(|y|^{k})$. It follows from the proof of Theorem \ref{teo1.1} that (\ref{1.5novoDois}) is formally conjugated to\begin{equation}\nonumber
		Y(y) = Ly+X_{K^3}(y)+\cdots+X_{K^{k-1}}(y)+X_{K^k}(y)+O(|y|^k),
	\end{equation}with $K^k\in{\mathcal{P}_{\tilde{\sigma}_1\tilde{\sigma}_2}^k(\mathbf{S}\rtimes\Gamma)}$. Since $\xi^k \in{\mathcal{P}}^k_{\sigma_1}(\Gamma)$, from Lemma \ref{mudsimpletica equivariante} we have that $\xi$ is $\Gamma$-equivariant. Therefore the new vector field preserves the symmetries (it is $\Gamma_{\sigma_2}$-equivariant) and the $\omega$-Hamiltonian condition of the original vector field. \cqd


\quad

Since $K^k\in{\mathcal{P}}^{k}_{\tilde{\sigma}_1\tilde{\sigma}_2}\left({\bf S}\rtimes\Gamma\right)\subset {\mathcal{P}}^{k}\left({\bf S}\right)$ for each $k\geq3$, the truncated equation\begin{equation}\label{1.14 novo novo}
	Y(y) = Ly+X_{K^3}(y)+\cdots+X_{K^{r}}(y)
\end{equation}is a normal form up to order $r-1$ of $X$. In fact, by (\ref{esta em Ck}) we have $X_{K^k}\in{\overrightarrow{{\mathcal{P}}}^{k-1}}(\mathbf{S})$. Moreover, $X_{K^k}\in{\overrightarrow{{\mathcal{P}}}_{\sigma_2}^{k-1}}(\Gamma)$, which implies by item 4 of Proposition \ref{propotrintadois} that $X_{K^k}\in {\overrightarrow{\mathcal{P}}_{\tilde{\sigma}_2}^{k-1}}(\mathbf{S}\rtimes\Gamma).$ So, Theorem \ref{fnormalnossa} gives us a normal form that preserves the symmetries of the original vector field. On the other hand, by Lemma \ref{LH2} the Hamiltonian function associated with the vector field (\ref{1.14 novo novo}) is given by $K = H^2+K^3+\cdots+K^{r}$, which is clearly a normal form of the Hamiltonian function $H$ associated with $X$, since $K^k\in\mathcal{P}^k(\mathbf{S})$. Hence we have the following result:

\begin{corollary}\label{fnormalnossa2}Under the conditions of the previous theorem, let $H:V\to\mathbb{R}$ be the Hamiltonian function associated with the vector field $X: V\to V$ such that $H(0)=0$ and $\nabla H(0)=0$. Let $H^2$ be the quadratic form associated with the Hessian of $H$. Then \begin{equation}\nonumber
		K = H^2+K^3+\cdots+K^{r},\end{equation}with $K^k\in {\mathcal{P}}^{k}_{\tilde{\sigma}_1\tilde{\sigma}_2}\left({\bf S}\rtimes\Gamma\right)$, is a normal form up to order $r$ of $H$ so that $X_K$ is a normal form up to order $r-1$ of the original vector field $X$ that preserves its symmetries.
\end{corollary}

Hence, in order to obtain a normal form of an $\omega$-Hamiltonian vector field under the action of a symmetry group $\Gamma$, it is enough to determine a normal form of the function $H$ according to Corollary \ref{fnormalnossa2}. In this case, we obtain the generators of degree $k$ of the module ${\mathcal{P}}_{\tilde{\sigma}_1\tilde{\sigma}_2}\left({\bf S}\rtimes\Gamma\right)$ over the ring ${\cal{P}}(\mathbf{S}\rtimes\Gamma)$, if $\sigma_1\sigma_2:\Gamma\to\mathbb{Z}_2$ is an epimorphism, or a Hilbert basis of the ring ${\mathcal{P}}_{\tilde{\sigma}_1\tilde{\sigma}_2}\left({\bf S}\rtimes\Gamma\right) = {\cal{P}}(\mathbf{S}\rtimes\Gamma)$, if $\sigma_1\sigma_2:\Gamma\to\mathbb{Z}_2$ is the trivial homomorphism. In practice, we exhibit a general form for the elements of ${\mathcal{P}}_{\tilde{\sigma}_1\tilde{\sigma}_2}\left({\bf S}\rtimes\Gamma\right)$, by using \cite[Theorem 3.1]{ABDM} and \cite[Theorem 3.2]{BM} if $\mathbf{S}\rtimes\Gamma$ is compact. When the group $\mathbf{S}\rtimes\Gamma$ is not compact, we can still use tools of invariant theory if ${\cal{P}}(\mathbf{S}\rtimes\Gamma)$ is finitely generated.

\quad

In the next two subsections, we apply the algebraic method presented in Theorem \ref{fnormalnossa} (and therefore in Corollary \ref{fnormalnossa2}) in order to obtain a normal form of two $\omega$-Hamiltonian vector fields.

\subsection{$\omega$-Hamiltonian $({\mathbb{D}}_4)_{\sigma}$-equivariant vector fields on $\mathbb{R}^2$}\label{sec:example}

Let $\Gamma =\mathbb{D}_4$ be the dihedral group generated by matrices \begin{equation}\nonumber
	R_{\frac{\pi}{2}} = \left[\begin{array}{cc}
		0&-1\\1&0
	\end{array}\right] \quad  \textrm{and} \quad   \kappa = \left[\begin{array}{cc}
		1&0\\0&-1
	\end{array}\right]
\end{equation}acting on an arbitrary symplectic space $({\mathbb{R}}^2,\omega)$. Consider \begin{equation}\nonumber
	[\omega] = \left[\begin{array}{cc}
		0&a_{12}\\-a_{12}&0
	\end{array}\right] \quad  \textrm{and} \quad L = \left[\begin{array}{cc}
		0&\lambda\\
		-\lambda&0
	\end{array}\right],
\end{equation} for some $a_{12}\neq0$ and $\lambda\in{\mathbb{R}}^*$. We want to determine a normal form of an $\omega$-Hamiltonian and $({\mathbb{D}}_4)_{\sigma}$-equivariant vector field $X:{\mathbb{R}}^2\to {\mathbb{R}}^2$ whose linearization at the origin is equal to $L$, where $\sigma:{\mathbb{D}}_4\to{\mathbb{Z}}_2$ is the epimorphism such that $\ker\sigma$ is generated by $R_{\frac{\pi}{2}}$. The vector field $X$ has Hamiltonian function $H: \mathbb{R}^2\to\mathbb{R}$ whose quadratic form is\begin{equation}\nonumber
	H^2(x_1,x_2) = \dfrac{\lambda a_{12}}{2}\left(x_1^2+x_2^2\right).
\end{equation} Clearly $R_{\frac{\pi}{2}}\in Sp_{\omega}(1,\mathbb{R})$ and $\kappa\in Sp_{\omega}^{-1}(1,\mathbb{R})$, that is, $\mathbb{D}_4$ is a subgroup of $\Omega_1$ as defined in (\ref{eq:omegan}). Moreover, $\ker\sigma\subset Sp_{\omega}(1,\mathbb{R})$. By Theorem \ref{acao subgrupo de Omega}, $\mathbb{D}_4$ acts $\sigma$-semisymplectically on $({\mathbb{R}}^2, \omega)$. This is a case where $\sigma_1 = \sigma_2: = \sigma$ and $X$ presents two types of symmetry of Table \ref{tabelatipos}: the elements of $\ker\sigma$ are SE and the elements of ${\mathbb{D}}_4\backslash\ker\sigma$ are AR.

In order to characterize the linear group $\mathbf{S}$, we observe that $e^{sL^T}$ is the rotation matrix
\begin{equation}\nonumber R_{\lambda s} = \left[
	\begin{array}{cc}
		\textrm{cos} \ \left(\lambda s\right) & -\textrm{sen} \ \left(\lambda s\right) \\
		\textrm{sen} \ \left(\lambda s\right) & \textrm{cos} \ \left(\lambda s\right) \\
	\end{array}
	\right],
\end{equation}whence $\left\{R_{\lambda s}: \ s\in\mathbb{R}\right\}$ is isomorphic to rotation group $SO(2)$. Therefore $\mathbf{S}=\overline{\left\{R_{\lambda s}: \ s\in\mathbb{R}\right\}}\simeq SO(2).$ In this way, by \cite[Example 3.1]{BRHT2021} and Theorem \ref{acao subgrupo de Omega} the action of \textbf{S} on $(\mathbb{R}^2,\omega)$ is symplectic, even if $a_{12}\neq1$ (non-symplectic system of coordinates).


By Theorem \ref{fnormalnossa}, it is enough to describe a generator set for the ring ${\mathcal{P}}_{\tilde{\sigma}_1\tilde{\sigma}_2}(SO(2)\rtimes \mathbb{D}_4) = {\mathcal{P}}(SO(2)\rtimes \mathbb{D}_4)$, since $\tilde{\sigma}_j(R_{\theta},\gamma) = \sigma_j(\gamma) = \sigma(\gamma)$, for $j=1,2$. By \cite[Theorem 3.2]{BM} and \cite[XII Examples 4.1]{Golu2}, a Hilbert basis of the ring ${\mathcal{P}}(SO(2)\rtimes \mathbb{D}_4)$ is $\{u(x_1,x_2) = x_1^2+x_2^2\}$. Hence, $H$ admits only normal forms of order $2r$ given by $$K(x_1,x_2) = H^2(x_1,x_2)+K^4(x_1,x_2)+\cdots+K^{2r}(x_1,x_2) = \dfrac{\lambda a_{12}}{2}(x_1^2+x_2^2)+\displaystyle\sum_{j=2}^{r}C_j(x_1^2+x_2^2)^j,$$ where $C_j\in\mathbb{R}$. Therefore a normal form up to order $2r-1$ of the $\omega$-Hamiltonian vector field $X$ is given by $X_K (x_1,x_2) = ([\omega]^{-1})^T\nabla K(x_1,x_2)$, that is,\begin{equation}\label{noplano}
	X_K (x_1,x_2) = \left(\lambda+\displaystyle\frac{2}{a_{12}}\sum_{j=2}^{r}jC_{j}(x_1^2+x_2^2)^{j-1}\right)(x_2, -x_1).
\end{equation}

From Theorem \ref{fnormalnossa}, the normal form (\ref{noplano}) preserves the two types of symmetry of the original vector field: symplectic equivariant and antisymplectic reversible. Moreover, by Theorem \ref{teo1.1} (without the action of a group $\Gamma$),  all $\omega$-Hamiltonian vector field on $({\mathbb{R}}^2,\omega)$ with linearization equal to $L$ has a normal form up to order $2r-1$ given by (\ref{noplano}). This happens because we can assume ${\mathcal{D}}^{k} = {\mathcal{P}^k}(SO(2))$, which coincides with ${\mathcal{P}^k}(SO(2)\rtimes\mathbb{D}_4)$, for all $k \geq 1$. In fact, by Proposition \ref{propotrintadois},\begin{equation}\nonumber
	{\mathcal{P}^k}(SO(2)\rtimes{\mathbb{D}}_4) = {\mathcal{P}^k}(SO(2))\cap{\mathcal{P}^k}({\mathbb{D}}_4) = {\mathcal{P}^k}(SO(2)),
\end{equation}since all $SO(2)$-invariant function generated by $u$ is also a ${\mathbb{D}}_4$-invariant function. Therefore, every $\omega$-Hamiltonian vector field on $\mathbb{R}^2$ with linearization $L$ is formally conjugated to an $\omega$-Hamiltonian vector field with SE and AR symmetries.

\subsection{$\omega$-Hamiltonian $({\mathbb{Z}}_2^{\tau}\times{\mathbb{Z}}_2^{\psi})_{\sigma_2}$-equivariant vector fields on $\mathbb{R}^4$}\label{sec:example2}

In Example \ref{exemploum}, we consider $\Gamma = {\mathbb{Z}}_2^{\tau}\times{\mathbb{Z}}_2^{\psi}$ acting $\sigma_1$-semisymplectically on $({\mathbb{R}}^4,\omega)$ and an $\omega$-Hamiltonian $\Gamma_{\sigma_2}$-equivariant vector field $X(x_1,x_2,x_3,x_4) = \lambda(x_2,-x_1,x_4,-x_3)$, for $\sigma_1, \sigma_2: \Gamma \to{\mathbb{Z}}_2$ such that $\ker\sigma_1 = \{(I_4,I_4), (I_4,\psi)\}$ and $\ker\sigma_2 = \{(I_4,I_4), (\tau,\psi)\}$.

Let $Y: {\mathbb{R}}^4 \to {\mathbb{R}}^4$ be an $\omega$-Hamiltonian and $\Gamma_{\sigma_2}$-equivariant vector field whose line\-arization at the origin is $L = dX_0$. In Table \ref{tabela4tipos}, we list the four types of symmetries of $Y$. By Example \ref{exemploum}, the quadratic form of the Hamiltonian function $H$ associated with $Y$ is given by $H^2(x_1,x_2,x_3,x_4) = -\lambda\left(x_1x_3+x_2x_4\right)$.

In order to determine a normal form for $Y$, we must describe the module ${\mathcal{P}}_{\tilde{\sigma}_1\tilde{\sigma}_2}({\bf S}\rtimes\Gamma)$ over the ring ${\mathcal{P}}({\bf S}\rtimes\Gamma)$. Since $\lambda i$ is the only algebraically independent eigenvalue of $L$ and this matrix has zero nilpotent part, by \cite[XVI Proposition 5.7]{Golu2} the group $\mathbf{S}$ is the circle group $S^1$ acting on $\mathbb{R}^4$ by 
$$\theta (x) = (x_1\cos \theta - x_2 \sin \theta, x_1 \sin \theta + x_2 \cos \theta, x_3\cos \theta - x_4\sin \theta, x_3 \sin \theta + x_4 \cos \theta),$$ for all $\theta \in S^1$ and $x = (x_1,x_2,x_3,x_4)\in \mathbb{R}^4$. A Hilbert basis for ${\mathcal{P}}(S^{\hspace{-0.05cm}1}\rtimes\Gamma)$ is given by $\{u_1,u_2,u_3^2,u_4^2\}$, where $u_1(x) = x_1^2 + x_2^2$, $u_2(x) = x_3^2 + x_4^2$, $u_3(x) = x_1x_3 + x_2x_4$ and $u_4(x) = x_2x_3 - x_1x_4$, and $u_3$ is a generator of the module ${\mathcal{P}}_{\tilde{\sigma}_1\tilde{\sigma}_2}(S^{\hspace{-0.05cm}1}\rtimes\Gamma)$ over ${\mathcal{P}}(S^{\hspace{-0.05cm}1}\rtimes\Gamma)$. For these computations, we use the tools presented in \cite[Theorem 3.2]{BM} and \cite[Theorem 3.1]{ABDM}, respectively. By Corollary \ref{fnormalnossa2}, a normal form up to order $2r$ of $H$ is given by $$K(x) = H^2(x)+\displaystyle\sum_{k=2}^{r} K^{2k}(x),$$ with $K^{2k}\in{\mathcal{{P}}_{\tilde{\sigma}_1\tilde{\sigma}_2}}(S^{\hspace{-0.05cm}1}\rtimes \Gamma)$ written as \begin{center}
$K^{2k}(x) = \displaystyle(x_1x_3+x_2x_4)\sum_{|\alpha|=k-1}C_{\alpha}(x_1^2+x_2^2)^{j_1}(x_3^2+x_4^2)^{j_2}(x_1x_3+x_2x_4)^{2j_3}(x_2x_3-x_1x_4)^{2j_4}$,
\end{center} for $\alpha = (j_1,j_2,2j_3,2j_4)$, $|\alpha|=j_1+j_2+2j_3+2j_4$, $C_{\alpha}\in{\mathbb{R}}$ and $k=2,\ldots,r$. Write 
\begin{align*}
F_1^k(x) &=\sum_{|\alpha| = k-1}2j_1C_{\alpha}u_1^{j_1-1}(x)u_2^{j_2}(x)u_3^{2j_3}(x)u_4^{2j_4}(x),\\
F_2^k(x)& = \sum_{|\alpha| = k-1}2j_2C_{\alpha}u_1^{j_1}(x)u_2^{j_2-1}(x)u_3^{2j_3}(x)u_4^{2j_4}(x),\\
F_3^k(x) &= \sum_{|\alpha| = k-1}(2j_3+1)C_{\alpha}u_1^{j_1}(x)u_2^{j_2}(x)u_3^{2(j_3-1)}(x)u_4^{2j_4}(x),\\
F_4^k(x)& =\sum_{|\alpha| = k-1}2j_4C_{\alpha}u_1^{j_1}(x)u_2^{j_2}(x)u_3^{2j_3}(x)u_4^{2(j_4-1)}(x)
\end{align*} for each $k=2,\ldots,r$. Thus \begin{align*}
\dfrac{\partial K^{2k}}{\partial x_1}(x) &= x_1u_3(x)F_1^k(x)+x_3u_3^2(x)F_3^k(x)-x_4u_3(x)u_4(x)F_4^k(x),\\
\dfrac{\partial K^{2k}}{\partial x_2}(x) &= x_2u_3(x)F_1^k(x)+x_4u_3^2(x)F_3^k(x)+x_3u_3(x)u_4(x)F_4^k(x),\\
\dfrac{\partial K^{2k}}{\partial x_3}(x) &= x_3u_3(x)F_2^k(x)+x_1u_3^2(x)F_3^k(x)+x_2u_3(x)u_4(x)F_4^k(x),\\ \dfrac{\partial K^{2k}}{\partial x_4}(x) &= x_4u_3(x)F_2^k(x)+x_2u_3^2(x)F_3^k(x)-x_1u_3(x)u_4(x)F_4^k(x).
\end{align*} Therefore, a normal form up to order $2r-1$ of the $\omega$-Hamiltonian vector field $Y$ is given by 
\begin{center}
$X_K(x) = Lx+[\omega]\nabla\left(\displaystyle\sum_{k=2}^{r}K^{2k}(x)\right) = (X_1(x),X_2(x),X_3(x),X_4(x)),$ 
\end{center} where $$
\begin{cases}
\displaystyle X_1(x) = \lambda x_2-x_4u_3(x)\sum_{k=2}^{r}F_2^k(x)-x_2u_3^2(x)\sum_{k=2}^{r}F_3^k(x)+x_1u_3(x)u_4(x)\sum_{k=2}^{r}F_4^k(x)\\
\displaystyle X_2(x) = -\lambda x_1+x_3u_3(x)\sum_{k=2}^{r}F_2^k(x)+x_1u_3^2(x)\sum_{k=2}^{r}F_3^k(x)+x_2u_3(x)u_4(x)\sum_{k=2}^{r}F_4^k(x)\\
\displaystyle X_3(x) = \lambda x_4-x_2u_3(x)\sum_{k=2}^{r}F_1^k(x)-x_4u_3^2(x)\sum_{k=2}^{r}F_3^k(x)-x_3u_3(x)u_4(x)\sum_{k=2}^{r}F_4^k(x)\\
\displaystyle X_4(x) = -\lambda x_3+x_1u_3(x)\sum_{k=2}^{r}F_1^k(x)+x_3u_3^2(x)\sum_{k=2}^{r}F_3^k(x)-x_4u_3(x)u_4(x)\sum_{k=2}^{r}F_4^k(x)
\end{cases}.$$

\section{Characterization of the complementary subspace of ${ad_{H^2}^k}(\mathcal{P}_{\sigma_1}^k(\Gamma))$ in $\mathcal{P}_{\sigma_1\sigma_2}^k(\Gamma)$}\label{adsubsecaopropriedades}
In this section, our goal is proving Theorem \ref{teo1.1nosso}, which characterizes a complementary subspace of ${ad_{H^2}^k}(\mathcal{P}_{\sigma_1}^k(\Gamma))$ in $\mathcal{P}_{\sigma_1\sigma_2}^k(\Gamma)$. We start by considering the action of $\Gamma$ on $(V,\omega)$, that induces for each $\gamma\in\Gamma$ the linear map $\rho_{\gamma}:V\to V$, $\rho_{\gamma}(x) = \gamma x$. Define the action ${}^*: \Gamma\times\mathcal{P}^k\to\mathcal{P}^k$ of $\Gamma$ on $\mathcal{P}^k$ as\begin{equation}
	\label{acaoPk}
	\gamma^*F = (\rho_{\gamma})^*F,
\end{equation}where $(\rho_{\gamma})^*F$ is the {pullback} of $F\in\mathcal{P}^k$ by $\rho_{\gamma}$, that is, $\gamma^*F(x) = F(\rho_{\gamma} (x)) = F(\gamma x),$ for all $\gamma\in\Gamma$ and $x\in V$. It follows from the linearity of the {pullback} of differential forms that the action (\ref{acaoPk}) is linear. Besides, if $L = dX_0$, then $L\gamma = \sigma_2(\gamma)\gamma L$ for all $\gamma\in\Gamma$, since $X\in{\overrightarrow{\mathcal{P}}_{\sigma_2}}(\Gamma)$.
\begin{lemma}\label{adequivariante}
	For each $k\geq3$, the operator $ad_{H^2}^k$ is $\Gamma_{\sigma_2}$-equivariant under the action (\ref{acaoPk}), that is, $ad_{H^2}^k(\gamma^*F) = \sigma_2(\gamma)\gamma^*ad_{H^2}^k(F)$, for all $\gamma\in\Gamma$ and $F\in\mathcal{P}^k$.
\end{lemma}
\Dem For each $x\in V$ and $F\in\mathcal{P}^k$, we have $ad_{H^2}^k(F)(x) = \{F,H^2\} (x) = dF_x(X_{H^2}(x)) = dF_xLx,$ where the last equality follows from Lemma \ref{LH2}. Thus, for each $\gamma\in\Gamma$ and $x\in V$, we have\begin{eqnarray}
		ad_{H^2}^k(\gamma^*F)(x) & = & d(\gamma^*F)_xLx = d(F(\rho_{\gamma}(x)))_xLx = dF_{\rho_{\gamma}(x)}d(\rho_{\gamma})_xLx \nonumber\\
		&=&dF_{\gamma x}\gamma Lx = \sigma_2(\gamma) dF_{\gamma x}L\gamma x = \sigma_2(\gamma)ad_{H^2}^k(F)(\gamma x)\nonumber\\
		& =& \sigma_2(\gamma)\gamma^*ad_{H^2}^k(F)(x).\nonumber
	\end{eqnarray} \cqd 

Next, we introduce some subgroups of $\Gamma$ that will be important in the proof of Lemmas \ref{lemaproj quase}-\ref{projecao de pS}. For group homomorphisms $\sigma_1,\sigma_2:\Gamma\to\mathbb{Z}_2$, fix $\delta_1, \delta_2, \delta_3\in\Gamma$ (if they there exist) such that
\begin{equation}
	\label{delta_}
	\sigma_j(\delta_i) =
	\begin{cases}
		\ \ 1,& \textrm{if} \ i=j\\
		-1,& \textrm{if} \ i\neq j
	\end{cases},
\end{equation}for $i=1,2,3$ and $j=1,2$. Denote
$$\Gamma_{++} = \ker\sigma_1\cap\ker\sigma_2,\quad \Gamma_{--} = (\Gamma\backslash\ker\sigma_1)\cap(\Gamma\backslash\ker\sigma_2),$$ $$\Gamma_{-+} = (\Gamma\backslash\ker\sigma_1)\cap\ker\sigma_2,\quad \Gamma_{+-} = \ker\sigma_1\cap(\Gamma\backslash\ker\sigma_2).$$ The signs ``$+$'' and ``$-$'' in the notations refer to the signs of the range of their elements by $\sigma_1$ and by $\sigma_2$, respectively. Then we have the following result:
\begin{proposition}\label{propo3coracao}
	Let $\sigma_1,\sigma_2:\Gamma\to\mathbb{Z}_2$ be group homomorphisms and $\delta_1, \delta_2, \delta_3\in\Gamma$ (if they there exist) satisfying (\ref{delta_}).\begin{enumerate}
		\item If $\sigma_1$ and $\sigma_2$ are distinct epimorphisms, then $\Gamma_{++}$ is a subgroup of $\Gamma$ of index $4$. Moreover, $\Gamma = \Gamma_{++} \ \dot{\cup} \ \Gamma_{+-} \ \dot{\cup} \ \Gamma_{-+} \ \dot{\cup} \ \Gamma_{--} = \Gamma_{++} \ \dot{\cup} \ \delta_1\Gamma_{++} \ \dot{\cup} \ \delta_2\Gamma_{++} \ \dot{\cup} \ \delta_3\Gamma_{++}.$
		\item If $\sigma_1$ is trivial and $\sigma_2$ is an epimorphism, then $\Gamma_{++} = \ker\sigma_2$ is a subgroup of $\Gamma$ of index $2$ and
		$\Gamma =  \Gamma_{++} \ \dot{\cup} \ \Gamma_{+-}  =  \Gamma_{++} \ \dot{\cup} \ \delta_1\Gamma_{++}.$
		\item If $\sigma_1$ is an epimorphism and $\sigma_2$ is trivial, then $\Gamma_{++} = \ker\sigma_1$ is a subgroup of $\Gamma$ of index $2$ and $\Gamma = \Gamma_{++} \ \dot{\cup} \ \Gamma_{-+} = \Gamma_{++} \ \dot{\cup} \ \delta_2\Gamma_{++}.$
		\item If $\sigma: = \sigma_1 = \sigma_2$ is an epimorphism, then $\Gamma_{++} = \ker\sigma$ is a subgroup of $\Gamma$ of index $2$ and we can write $\Gamma  = \Gamma_{++} \ \dot{\cup} \ \Gamma_{--}  =  \Gamma_{++} \ \dot{\cup} \ \delta_3\Gamma_{++}.$
		\item If $\sigma_1 = \sigma_2$ is the trivial homomorphism, then $\Gamma = \Gamma_{++}$. 
	\end{enumerate}
\end{proposition}
\Dem We only prove the first item, since the other ones follow similarly.
	If $\sigma_1$ and $\sigma_2$ are distinct group epimorphisms, then by the proof of Proposition \ref{proptipos} there exist $\delta_1$, $\delta_2$ and $\delta_3 = \delta_1\delta_2$ in $\Gamma$ satisfying (\ref{delta_}), that is, $\delta_1\in \Gamma_{+-}$, $\delta_2\in\Gamma_{-+}$ and $\delta_3\in\Gamma_{--}$. We now prove that $\Gamma_{++}$ is a subgroup of $\Gamma$ of index $[\Gamma:\Gamma_{++}] = 4$. In fact, consider the restricted map $\sigma_1\mid_{\ker\sigma_2}:\ker\sigma_2\to\mathbb{Z}_2$. Clearly $\ker(\sigma_1\mid_{\ker\sigma_2}) \subset\ker\sigma_1\cap\ker\sigma_2= \Gamma_{++}$. On the other hand, if $\gamma\in\Gamma_{++}$ then $\gamma\in\ker\sigma_1\cap\ker\sigma_2$. Thus $\sigma_1\mid_{\ker\sigma_2}(\gamma)$ is well defined and\begin{equation}\nonumber
		\sigma_1\mid_{\ker\sigma_2}(\gamma) = \sigma_1(\gamma)  = 1,
	\end{equation} which implies that $ \Gamma_{++}\subset\ker(\sigma_1\mid_{\ker\sigma_2})$. Therefore, $\Gamma_{++}= \ker\left(\sigma_1\mid_{\ker\sigma_2}\right)$ is a subgroup of $\ker\sigma_2$ of index $2$. Hence\begin{equation}\nonumber
		[\Gamma:\Gamma_{++}] = [\Gamma:\ker\sigma_2]\cdot [\ker\sigma_2:\Gamma_{++}] = 4.\end{equation}
	In this case, the left-cosets of $\Gamma_{++}$ in $\Gamma$ are $\Gamma_{++}$, $ \Gamma_{+-}=\delta_1\Gamma_{++}$, $\Gamma_{-+}=\delta_2\Gamma_{++}$ and $\Gamma_{--} = \delta_3\Gamma_{++}$, which completes the proof.
\cqd

An important tool in representation theory of groups is the Haar integral, an invariant form of integration by translation of elements of a Lie group. Based on the approach presented in \cite[XVI Theorem 5.8; Theorem 5.9]{Golu2} and \cite[Theorem 4.7]{Bap}, we define a projection of the ring $\mathcal{P}^k(\mathbf{S})$ of all \textbf{S}-invariant polynomial functions onto the module ${\mathcal{P}}_{\tilde{\sigma}_1\tilde{\sigma}_2}^{k}(\mathbf{S}\rtimes\Gamma)$ (Lemma \ref{projecao de pS}), where the action of $\mathbf{S}$ on $(V,\omega)$ is symplectic and $\tilde{\sigma}_j:\mathbf{S}\rtimes\Gamma\to\mathbb{Z}_2$ is defined as in (\ref{estrelaCorrecPatMuitasvezes}).
 

Suppose $\Gamma_{++} = \ker\sigma_1\cap\ker\sigma_2$ is a linear Lie group. In this case, $\Gamma_{++}$ is closed in $\Gamma$ and, therefore, compact. This hypothesis is important in the computation of the normalized Haar integral over $\Gamma_{++}$. 

\begin{definition}
	Consider the action ${}^*: \Gamma\times\mathcal{P}^k\to\mathcal{P}^k$ defined in (\ref{acaoPk}). For each $k\geq3$, define the maps $\overline{\pi}, \pi:{\mathcal{P}^k}\to{\mathcal{P}^k}$ as
	$$		\overline{\pi}(F)=\displaystyle\frac{1}{[\Gamma:\Gamma_{++}]}\left(\int_{\tau\in\Gamma_{++}}\tau^*F+\int_{\tau\in\Gamma_{++}}(\delta_1\tau)^*F - \int_{\tau\in\Gamma_{++}}(\delta_2\tau)^*F-\int_{\tau\in\Gamma_{++}}(\delta_3\tau)^*F\right)
			$$ and $$
			\pi(F)=\displaystyle\frac{1}{[\Gamma:\Gamma_{++}]}\left(\int_{\tau\in\Gamma_{++}}\tau^*F-\int_{\tau\in\Gamma_{++}}(\delta_1\tau)^*F - \int_{\tau\in\Gamma_{++}}(\delta_2\tau)^*F+\int_{\tau\in\Gamma_{++}}(\delta_3\tau)^*F\right),$$ where $\displaystyle\int_{\tau\in\Gamma_{++}}$ is the normalized Haar integral over $\Gamma_{++}$. We assume $\displaystyle\int_{\tau\in\Gamma_{++}}(\delta_i\tau)^*F = 0$ if $\delta_i$ does not exist, $i=1,2,3$. \end{definition}
\begin{lemma}\label{lemaproj quase}For each $k\geq3$, the map $\overline{\pi}:{\mathcal{P}^k}\to{\mathcal{P}_{\sigma_1}^k}(\Gamma)$ is an idempotent linear projection, that is, $\overline{\pi}^2 = \overline{\pi}$.
\end{lemma}
\Dem Since the {pullback} and the Haar integral are linear, we have that $\overline{\pi}$ is linear. For the same reason, if the degree of $F$ is $k$, then the degree of $\overline{\pi}(F)$ is also $k$. We want to show that $\overline{\pi}(F) \in{\mathcal{P}_{\sigma_1}^k}(\Gamma)$, that is, $\gamma^*\overline{\pi}(F) = \sigma_1(\gamma)\overline{\pi}(F)$, for any $\gamma\in\Gamma$ and $F\in\mathcal{P}^k$. Based on Proposition \ref{propo3coracao}, we prove the case where $\sigma_1$ and $\sigma_2$ are distinct epimorphisms, that is, $[\Gamma:\Gamma_{++}] = 4$. The other cases follow similarly. 
	
	Given $\gamma\in\Gamma$,
	$$\begin{array}{rcl}
			\gamma^*\overline{\pi}(F) &= &\gamma^*\left(\displaystyle\frac{1}{4}\left(\int_{\tau\in\Gamma_{++}}\tau^*F+\int_{\tau\in\Gamma_{++}}(\delta_1\tau)^*F - \int_{\tau\in\Gamma_{++}}(\delta_2\tau)^*F-\int_{\tau\in\Gamma_{++}}(\delta_3\tau)^*F\right)\right)\vspace{0.15cm}\\
			
			& = & \displaystyle\frac{1}{4}\left(\int_{\tau\in\Gamma_{++}}\gamma^*\tau^*F+\int_{\tau\in\Gamma_{++}}\gamma^*(\delta_1\tau)^*F - \int_{\tau\in\Gamma_{++}}\gamma^*(\delta_2\tau)^*F -\int_{\tau\in\Gamma_{++}}\gamma^*(\delta_3\tau)^*F \right)\vspace{0.15cm}\\
			
			& = & \displaystyle\frac{1}{4}\left(\int_{\tau\in\Gamma_{++}}(\tau\gamma)^*F +\int_{\tau\in\Gamma_{++}}(\delta_1\tau\gamma)^*F - \int_{\tau\in\Gamma_{++}}(\delta_2\tau\gamma)^*F -\int_{\tau\in\Gamma_{++}}(\delta_3\tau\gamma)^*F \right),
		\end{array} $$
\noindent where in the second equality we use the continuity of the {pullback} and in the third equality we use that $(\tau\gamma)^*F = \gamma^*\tau^*F$ for all $\tau, \gamma\in\Gamma$ and $F\in\mathcal{P}^k$.
	
	If $\gamma\in\Gamma_{++}$, then $\tau\gamma\in\Gamma_{++}$. Since the Haar integral is invariant by elements of $\Gamma_{++}$ we have that $\gamma^*\overline{\pi}(F) = \overline{\pi} (F) = \sigma_1(\gamma)\overline{\pi}(F)$.
	
	If $\gamma\in\Gamma_{+-} = \delta_1\Gamma_{++}$, then there exists $\bar{\gamma}\in\Gamma_{++}$ such that $\gamma = \delta_1\bar{\gamma}$. Note that $\tau\delta_1\in\Gamma_{+-}= \delta_1\Gamma_{++}$, $\delta_3\delta_1\in\Gamma_{-+}= \delta_2\Gamma_{++}$ and $\delta_2\delta_1\in\Gamma_{--}= \delta_3\Gamma_{++}$, because $\sigma_j(\tau\delta_1) = \sigma_j(\delta_1)$, $\sigma_j(\delta_3\delta_1) = - \sigma_j(\delta_1) = \sigma_j(\delta_2)$ and $\sigma_j(\delta_2\delta_1) = -1=\sigma_j(\delta_3)$, for all $j=1,2$. Hence we conclude that there exist $\bar{\gamma}, \bar{\tau}, \bar{\delta}_2, \bar{\delta}_3\in\Gamma_{++}$ such that $\gamma = \delta_1\bar{\gamma},$ $\tau\delta_1 = \delta_1\bar{\tau},$ $\delta_3\delta_1 = \delta_2\bar{\delta}_2$ and $\delta_2\delta_1 = \delta_3\bar{\delta}_3.$ In this case, by the invariance of the Haar integral by $\Gamma_{++}$ we have
$$\begin{array}{rcl}
			\gamma^*\overline{\pi}(F) & = & \displaystyle\frac{1}{4}\left(\int_{\tau\in\Gamma_{++}}(\tau\delta_1\bar{\gamma})^*F +\int_{\tau\in\Gamma_{++}}(\delta_1\tau\delta_1\bar{\gamma})^*F - \int_{\tau\in\Gamma_{++}}(\delta_2\tau\delta_1\bar{\gamma})^*F-\int_{\tau\in\Gamma_{++}}(\delta_3\tau\delta_1\bar{\gamma})^*F\right)\vspace{0.15cm}\\
			& = & \displaystyle\frac{1}{4}\left(\int_{\bar{\tau}\in\Gamma_{++}}(\delta_1\bar{\tau}\bar{\gamma})^*F+\int_{\bar{\tau}\in\Gamma_{++}}(\delta_1^2\bar{\tau}\bar{\gamma})^*F - \int_{\bar{\tau}\in\Gamma_{++}}(\delta_2\delta_1\bar{\tau}\bar{\gamma})^*F -\int_{\bar{\tau}\in\Gamma_{++}}(\delta_3\delta_1\bar{\tau}\bar{\gamma})^*F \right)\vspace{0.15cm}\\
			& = & \displaystyle\frac{1}{4}\left(\int_{\bar{\tau}\in\Gamma_{++}}(\delta_1\bar{\tau})^*F +\int_{\bar{\tau}\in\Gamma_{++}}\bar{\tau}^*F - \int_{\bar{\tau}\in\Gamma_{++}}(\delta_3\bar{\delta}_3\bar{\tau})^*F-\int_{\bar{\tau}\in\Gamma_{++}}(\delta_2\bar{\delta}_2\bar{\tau})^*F \right)\vspace{0.15cm}\\
			& = & \displaystyle\frac{1}{4}\left(\int_{\bar{\tau}\in\Gamma_{++}}(\delta_1\bar{\tau})^*F +\int_{\bar{\tau}\in\Gamma_{++}}\bar{\tau}^*F - \int_{\bar{\tau}\in\Gamma_{++}}(\delta_3\bar{\tau})^*F -\int_{\bar{\tau}\in\Gamma_{++}}(\delta_2\bar{\tau})^*F \right)\vspace{0.15cm}\\
			&= &\overline{\pi}(F) = \sigma_1(\gamma)\overline{\pi}(F).
		\end{array} $$
	
	If $\gamma\in\Gamma_{-+} = \delta_2\Gamma_{++}$, since $\tau\delta_2\in\Gamma_{-+}= \delta_2\Gamma_{++}$, $\delta_3\delta_2\in\Gamma_{+-}= \delta_1\Gamma_{++}$ and $\delta_1\delta_2\in\Gamma_{--}= \delta_3\Gamma_{++}$, then there exist  $\bar{\gamma}, \bar{\tau}, \bar{\delta}_1,\bar{\delta}_3\in\Gamma_{++}$ such that $\gamma = \delta_2\bar{\gamma},$ $\tau\delta_2 = \delta_2\bar{\tau},$ $\delta_3\delta_2 = \delta_1\bar{\delta}_1$ and $\delta_1\delta_2 = \delta_3\bar{\delta}_3.$ Again, by the invariance of the Haar integral by $\Gamma_{++}$, we have similarly to the previous case that
$$\begin{array}{rcl}
			\gamma^*\overline{\pi}(F) & = & \displaystyle\frac{1}{4}\left(\int_{\tau\in\Gamma_{++}}(\tau\delta_2\bar{\gamma})^*F+\int_{\tau\in\Gamma_{++}}(\delta_1\tau\delta_2\bar{\gamma})^*F - \int_{\tau\in\Gamma_{++}}(\delta_2\tau\delta_2\bar{\gamma})^*F-\int_{\tau\in\Gamma_{++}}(\delta_3\tau\delta_2\bar{\gamma})^*F \right)\vspace{0.15cm}\\
			& = & \displaystyle\frac{1}{4}\left(\int_{\bar{\tau}\in\Gamma_{++}}(\delta_2\bar{\tau})^*F +\int_{\bar{\tau}\in\Gamma_{++}}(\delta_3\bar{\tau})^*F - \int_{\bar{\tau}\in\Gamma_{++}}\bar{\tau}^*F -\int_{\bar{\tau}\in\Gamma_{++}}(\delta_1\bar{\tau})^*F \right) \vspace{0.15cm}\\
			&= &  - \overline{\pi}(F) = \sigma_1(\gamma)\overline{\pi}(F).
		\end{array} $$
		
If $\gamma\in\Gamma_{--} = \delta_3\Gamma_{++}$, since $\tau\delta_3\in\Gamma_{--}= \delta_3\Gamma_{++}$, $\delta_2\delta_3\in\Gamma_{+-}= \delta_1\Gamma_{++}$ and $\delta_1\delta_3\in\Gamma_{-+}= \delta_2\Gamma_{++}$, then there exist $\bar{\gamma}, \bar{\tau}, \bar{\delta}_1, \bar{\delta}_2\in\Gamma_{++}$ such that $\gamma = \delta_3\bar{\gamma},$ $\tau\delta_3 = \delta_3\bar{\tau},$ $\delta_2\delta_3 = \delta_1\bar{\delta}_1$ and $\delta_1\delta_3 = \delta_2\bar{\delta}.$ Hence we obtain
$$\begin{array}{rcl}
			\gamma^*\overline{\pi}(F) & = & \displaystyle\frac{1}{4}\left(\int_{\tau\in\Gamma_{++}}(\tau\delta_3\bar{\gamma})^*F +\int_{\tau\in\Gamma_{++}}(\delta_1\tau\delta_3\bar{\gamma})^*F - \int_{\tau\in\Gamma_{++}}(\delta_2\tau\delta_3\bar{\gamma})^*F -\int_{\tau\in\Gamma_{++}}(\delta_3\tau\delta_3\bar{\gamma})^*F \right)\vspace{0.15cm}\\
			& = & \displaystyle\frac{1}{4}\left(\int_{\bar{\tau}\in\Gamma_{++}}(\delta_3\bar{\tau})^*F +\int_{\bar{\tau}\in\Gamma_{++}}(\delta_2\bar{\tau})^*F - \int_{\bar{\tau}\in\Gamma_{++}}(\delta_1\bar{\tau})^*F -\int_{\bar{\tau}\in\Gamma_{++}}\bar{\tau}^*F\right)\vspace{0.15cm}\\	
			&= &- \overline{\pi}(F)  = \sigma_1(\gamma)\overline{\pi}(F).
		\end{array} $$ Therefore, in each case $\overline{\pi}(F)\in\mathcal{P}_{\sigma_1}^k(\Gamma)$. Since $F\in\mathcal{P}^k$ is arbitrary, we conclude that $\overline{\pi}(\mathcal{P}^k)\subset\mathcal{P}_{\sigma_1}^k(\Gamma)$. Reciprocally, given a function $F\in\mathcal{P}_{\sigma_1}^k(\Gamma)$, we have by (\ref{acaoPk}) that $\gamma^*F(x) = F(\gamma x) = \sigma_1(\gamma) F(x)$ for all $\gamma\in\Gamma$ and $x\in V$. Then
	$$\begin{array}{rcl}
			\overline{\pi}(F) & = &\displaystyle\frac{1}{4}\left(\int_{\tau\in\Gamma_{++}}\tau^*F +\int_{\tau\in\Gamma_{++}}(\delta_1\tau)^*F  -\int_{\tau\in\Gamma_{++}}(\delta_2\tau)^*F  -\int_{\tau\in\Gamma_{++}}(\delta_3\tau)^*F \right)\vspace{0.15cm}\\
			
			& = &\displaystyle\frac{1}{4}\left(\int_{\tau\in\Gamma_{++}}\sigma_1(\tau) F +\int_{\tau\in\Gamma_{++}}\sigma_1(\delta_1\tau)F  -\int_{\tau\in\Gamma_{++}}\sigma_1(\delta_2\tau)F  -\int_{\tau\in\Gamma_{++}}\sigma_1(\delta_3\tau)F \right)\vspace{0.15cm}\\
			
			& = &\displaystyle\frac{1}{4}\left(\int_{\tau\in\Gamma_{++}}F +\int_{\tau\in\Gamma_{++}}F -\int_{\tau\in\Gamma_{++}} (- F) -\int_{\tau\in\Gamma_{++}}(- F)\right)\vspace{0.15cm}\\
			
			& = & F\displaystyle \int_{\tau\in\Gamma_{++}}1 = F,
		\end{array}$$ \noindent where in the third equality we have the Haar integral of the constant function $f(\gamma) = F$. Thus $F\in\overline{\pi}(\mathcal{P}^k)$ and we conclude that $\overline{\pi}(\mathcal{P}^k) = \mathcal{P}_{\sigma_1}^k(\Gamma)$. Moreover $\overline{\pi}\mid_{\mathcal{P}_{\sigma_1}^k(\Gamma)} = {\rm Id}$. Therefore, $\overline{\pi}$ is an idempotent projection.\cqd

The proof of the next lemma is similar to the previous one.

\begin{lemma}\label{lemaproj principal}For each $k\geq3$, the map $\pi:{\mathcal{P}^k}\to{\mathcal{P}_{\sigma_1\sigma_2}^k}(\Gamma)$ is an idempotent linear projection.
\end{lemma}
 
\begin{lemma}\label{projecao de IMad}For each $k\geq3$, we have $\pi(ad_{H^2}^k(\mathcal{P}^k)) = ad_{H^2}^k(\overline{\pi}(\mathcal{P}^k)) = ad_{H^2}^k(\mathcal{P}_{\sigma_1}^k(\Gamma)).$
\end{lemma}
\Dem As before, we prove the case where $\sigma_1$ and $\sigma_2$ are distinct epimorphisms. The other cases are similar. Moreover, to simplify the notation, we will omit $\tau\in\Gamma_{++}$ in the Haar integral over $\Gamma_{++}$.

Given $F\in\mathcal{P}^k$, the linearity of the operator $ad_{H^2}^k$ ensures that\begin{equation}\nonumber
		\int ad_{H^2}^k(\tau^*F)  = ad_{H^2}^k\left(\int \tau^*F\right).
	\end{equation}Furthermore, by Lemma \ref{adequivariante}, we have $\gamma^*ad_{H^2}^k(F) = \sigma_2(\gamma)ad_{H^2}^k(\gamma^*F)$ for all $\gamma\in\Gamma$. Since $\sigma_1$ and $\sigma_2$ are distinct, there exist $\delta_1, \delta_2, \delta_3\in\Gamma$ satisfying (\ref{delta_}). For all $\tau\in\Gamma_{++}$, we have $\sigma_2(\delta_1\tau) = \sigma_2(\delta_3\tau) = -\sigma_2(\delta_2\tau) = -\sigma_2(\tau) = -1$. Moreover, $[\Gamma:\Gamma_{++}] = 4$. So, for all $F\in\mathcal{P}^k$,
	\begin{align*}
	\pi(ad_{H^2}^k(F)) & = \frac{1}{4}\left(\int \tau^*ad_{H^2}^k(F)-\int (\delta_1\tau)^*ad_{H^2}^k(F) - \displaystyle\int (\delta_2\tau)^*ad_{H^2}^k(F)+\int(\delta_3\tau)^*ad_{H^2}^k(F)\right)\\
		&=\displaystyle\frac{1}{4}\left(\int \sigma_2(\tau)ad_{H^2}^k(\tau^*F)-\int \sigma_2(\delta_1\tau)ad_{H^2}^k((\delta_1\tau)^*F) \right.\\
		&\left.\displaystyle \quad - \int \sigma_2(\delta_2\tau) ad_{H^2}^k((\delta_2\tau)^*F)+\int \sigma_2(\delta_3\tau)ad_{H^2}^k((\delta_3\tau)^*F)\right)\\
		&=\displaystyle\frac{1}{4}\left(\int ad_{H^2}^k(\tau^*F)+\int ad_{H^2}^k((\delta_1\tau)^*F) -\int ad_{H^2}^k((\delta_2\tau)^*F)-\int ad_{H^2}^k((\delta_3\tau)^*F)\right)\\
		&=ad_{H^2}^k(\overline{\pi}(F)). 
	\end{align*}
	Thus we obtain $\pi(ad_{H^2}^k(\mathcal{P}^k)) = ad_{H^2}^k(\overline{\pi}(\mathcal{P}^k))$. By Lemma \ref{lemaproj quase}, we have $\pi(ad_{H^2}^k(\mathcal{P}^k)) = ad_{H^2}^k(\mathcal{P}_{\sigma_1}^k(\Gamma))$, since $\overline{\pi}$ is surjective. \cqd
	
For the next result, we consider again $\tilde{\sigma}_j:\mathbf{S}\rtimes\Gamma\to\mathbb{Z}_2$ as in (\ref{estrelaCorrecPatMuitasvezes}), with $j=1,2$. We define in (\ref{6barra}) a group homomorphism on semidirect product $\Gamma_1\rtimes\Gamma_2$. In our case, $\tilde{\sigma}_j(e^{sL^T},\gamma) = \beta_{j1}(e^{sL^T})\beta_{j2}(\gamma)$, where $\beta_{j1}:\mathbf{S}\to\mathbb{Z}_2$ is the trivial homomorphism and $\beta_{j2}:\Gamma\to\mathbb{Z}_2$ coincides with $\sigma_j$. Furthermore,\begin{eqnarray}
	(\tilde{\sigma}_1\tilde{\sigma}_2)(e^{sL^T},\gamma)=\tilde{\sigma}_1(e^{sL^T},\gamma)\tilde{\sigma}_2(e^{sL^T},\gamma) = (\sigma_1\sigma_2)(\gamma),\nonumber
\end{eqnarray}for all $s\in\mathbb{R}$ and $\gamma\in\Gamma$.
\begin{lemma}\label{projecao de pS}For each $k\geq3$, we have $\pi(\mathcal{P}^k(\mathbf{S})) = \mathcal{P}_{\tilde{\sigma}_1\tilde{\sigma}_2}^k(\mathbf{S}\rtimes\Gamma)$.
\end{lemma}
\Dem By item 3 of Proposition \ref{propotrintadois}, \begin{equation}\label{item(iii)}
		{{\mathcal{P}}_{\tilde{\sigma}_1\tilde{\sigma}_2}}(\mathbf{S}\rtimes{\Gamma}) = {{\mathcal{P}}}(\mathbf{S}) \cap {{\mathcal{P}}_{\sigma_1\sigma_2}}({\Gamma}).
	\end{equation} Given $F\in\mathcal{P}_{\tilde{\sigma}_1\tilde{\sigma}_2}^k(\mathbf{S}\rtimes\Gamma)$, it follows that $F\in\mathcal{P}^k(\mathbf{S})$ and $F\in\mathcal{P}_{{\sigma}_1{\sigma}_2}^k(\Gamma)$. From proof of Lemma \ref{lemaproj principal} we have $F=\pi(F)\in\pi\left(\mathcal{P}^k(\mathbf{S})\right)$. Thus $\mathcal{P}_{\tilde{\sigma}_1\tilde{\sigma}_2}^k(\mathbf{S}\rtimes\Gamma)\subset\pi(\mathcal{P}^k(\mathbf{S}))$.
	
	On the other hand, given $F\in\mathcal{P}^k(\mathbf{S})$, we want to show that $\pi(F)\in\mathcal{P}_{\tilde{\sigma}_1\tilde{\sigma}_2}^k(\mathbf{S}\rtimes\Gamma)$, that is,\begin{equation}\nonumber\pi(F)((e^{sL^T},\gamma) x) = (\tilde{\sigma}_1\tilde{\sigma}_2)(e^{sL^T},\gamma) \pi(F)(x)= (\sigma_1\sigma_2)(\gamma)\pi(F)(x),
	\end{equation}for all $\gamma\in\Gamma, \ x\in V$. By (\ref{doisasterisco}), it is enough to show that
	\begin{equation}\label{umbola}\left(\gamma e^{\sigma_2(\gamma)sL^T}\right)^*\pi(F) = (\sigma_1\sigma_2)(\gamma)\pi(F),
	\end{equation}for all $\gamma\in\Gamma$ and $s\in\mathbb{R}$, where ${}^*$ is the action of $\mathbf{S}\rtimes\Gamma$ on $\mathcal{P}^k$ defined as in (\ref{acaoPk}). By (\ref{doisasterisco}), $e^{\sigma_2(\gamma)sL^T} = \gamma e^{sL^T}\gamma^{-1}$ and, in particular, for all $\gamma\in\ker\sigma_2$ we have $\gamma e^{sL^T} = e^{sL^T}\gamma$. Thus if $\gamma\in\Gamma$, $s\in\mathbb{R}$ and $x\in V$, then
	\begin{align}\label{++++2}
		(\gamma e^{\sigma_2(\gamma)sL^T})^* F(x)& =  F(\gamma e^{\sigma_2(\gamma)sL^T}x)  =  F(\gamma^2 e^{sL^T}\gamma^{-1} x)\nonumber \\
		& =  F(e^{sL^T}\gamma^2\gamma^{-1} x) =  F(\gamma x)\nonumber\\
		&   =   \gamma^*F(x),
	\end{align} where the third and fourth equalities follow since $\gamma^2\in\ker\sigma_2$ and $F\in\mathcal{P}^k(\mathbf{S})$, respectively. Thus, for all $\bar{\gamma}\in\Gamma$ and $x\in V$, we have
	\begin{align*}
			(\bar{\gamma}\gamma e^{\sigma_2(\gamma)sL^T})^*F(x)& = (\gamma e^{\sigma_2(\gamma)sL^T})^*((\bar{\gamma})^*F)(x) = (\gamma e^{\sigma_2(\gamma)sL^T})^*F(\bar{\gamma} x)\\
			& = \gamma^*F(\bar{\gamma} x) =  \gamma^*(\bar{\gamma})^*F(x) = (\bar{\gamma}\gamma)^*F(x),
		\end{align*}where the third equality follows from (\ref{++++2}). 
	
	As before, we present a proof only when $\sigma_1$ and $\sigma_2$ are distinct epimorphisms and omit $\tau\in\Gamma_{++}$ in the Haar integral. In this case, for all $\gamma\in\Gamma$ and $s\in\mathbb{R}$, we have\begin{align*}
		\left(\gamma e^{\sigma_2(\gamma)sL^T}\right)^*\pi(F) &=\left(\gamma e^{\sigma_2(\gamma)sL^T}\right)^*\left(\frac{1}{4} \left( \int \tau^*F - \int (\delta_1\tau)^*F - \int (\delta_2\tau)^*F + \int (\delta_3\tau)^*F\right) \right)\\ 
		&=\frac{1}{4}\left(\int (\tau\gamma e^{\sigma_2(\gamma)sL^T})^*F-\int (\delta_1\tau\gamma e^{\sigma_2(\gamma)sL^T})^*F \right.\\ &  \left.  \quad - \int (\delta_2\tau\gamma e^{\sigma_2(\gamma)sL^T})^*F + \int (\delta_3\tau\gamma e^{\sigma_2(\gamma)sL^T})^*F\right) \\
		&=\frac{1}{4}\left(\int (\tau\gamma)^*F-\int (\delta_1\tau\gamma)^*F - \int (\delta_2\tau\gamma )^*F+\int (\delta_3\tau\gamma)^*F\right)\\
		&=\gamma^*\pi(F) = (\sigma_1\sigma_2)(\gamma)\pi(F),
	\end{align*}where the last equality follows from Lemma \ref{lemaproj principal}. Therefore, the equation (\ref{umbola}) is valid. \cqd 

We are now able to present a proof for Theorem \ref{teo1.1nosso}.\\

\noindent {\sc Proof of Theorem \ref{teo1.1nosso}:} Since $\pi:\mathcal{P}^k\to\mathcal{P}^k$ is  linear, we can apply it to both sides of decomposition (\ref{complemento elphick}) in order to get\begin{equation}\nonumber
		\pi(\mathcal{P}^k) = \pi(ad_{H^2}^k(\mathcal{P}^k))+\pi(\mathcal{P}^k(\mathbf{S})).
	\end{equation}From Lemmas \ref{lemaproj principal}, \ref{projecao de IMad} and \ref{projecao de pS}, respectively, we have $${\mathcal{P}}^{k}_{\sigma_1\sigma_2}(\Gamma) = ad_{H^2}^{k}({\mathcal{P}}^{k}_{\sigma_1}(\Gamma))+{\mathcal{P}}^{k}_{\tilde{\sigma}_1\tilde{\sigma}_2}\left({\bf S}\rtimes\Gamma\right),$$ with $\tilde{\sigma}_j$ as defined in (\ref{estrelaCorrecPatMuitasvezes}), $j=1,2$. From (\ref{item(iii)}) we know that ${\mathcal{P}}^{k}_{\tilde{\sigma}_1\tilde{\sigma}_2}\left({\bf S}\rtimes\Gamma\right) =  {\mathcal{P}}^{k}(\mathbf{S}) \cap {\mathcal{P}}_{\sigma_1\sigma_2}^{k}(\Gamma)$ and by decomposition (\ref{complemento elphick}) we have that $ad_{H^2}^{k}({\mathcal{P}}^{k})\cap {\mathcal{P}}^{k}(\mathbf{S})=\left\{0\right\}$, which completes the proof. \cqd

\noindent {\bf Acknowledgments.} Eralcilene Moreira Terezio was financed in part by the Coordena\c{c}\~ao de Aperfei\c{c}oamento de Pessoal de N\'ivel Superior - Brasil (CAPES) - Finance 001. The research of Patrícia Hernandes Baptistelli was supported by FAPESP, grant 2019/21181-0.

\end{document}